\documentclass{article}
\usepackage{amsmath}
\usepackage{multicol,color}
\usepackage{float}
\usepackage{soul}
\usepackage{graphicx}
\usepackage{amssymb}
\usepackage{theorem}
\usepackage{fancyhdr}
\usepackage{tikz}
\usepackage{enumerate}
\usepackage[margin=1.1in]{geometry}
\def\disp{\displaystyle}
\def\ve{\varepsilon}
\def\dd{\delta}

\def\lm{\lambda}
\def\O{\Omega}
\def\Tilde{\widetilde}

\def\oa{\bar a}

\def\os{\bar s}
\def\ox{\bar{x}}
\def\oy{\bar{y}}
\def\oz{\bar{z}}

\def\ou{\bar{u}}

\def\os{\bar{s}}

\def\cone{\hbox{}}
\def\span{\hbox{\rm span}}
\def\gph{\hbox{}}

\def\gg{\gamma}
\def\dn{\downarrow}

\def\tto{\rightrightarrows}
\def\st{\stackrel}

\def\Hat{\widehat}

\def\Tilde{\widetilde}
\def\Bar{\overline}
\def\ra{\rangle}
\def\la{\langle}
\def\ve{\varepsilon}

\def\h{\hfill\Box}
\def\R{\mathbb{R}}
\def\N{\mathbb{N}}

\def\co{\mbox{\rm co}\,}
\def\gph{\mbox{\rm gph}\,}

\def\dom{\mbox{\rm dom}\,}

\def\dist{\mbox{\rm dist}}

\def\cone{\mbox{\rm cone}\,}

\def\var{\mbox{\rm var}\,}
\def\dn{\downarrow}
\def\O{\Omega}
\def\bP{({\bar P}^\t)}

\def\ph{\varphi}
\def\emp{\emptyset}
\def\st{\stackrel}
\def\oR{\Bar{\R}}
\def\lm{\lambda}

\def\gg{\gamma}

\def\dd{\delta}

\def\al{\alpha}
\def\vt{\vartheta}
\def\t{\tau}

\def\N{I\!\!N}
\def\th{\theta}
\def\vTh{\vartheta}
\newtheorem{theorem}{Theorem}[section]

\newtheorem{corollary}[theorem]{Corollary}
\newtheorem{proposition}[theorem]{Proposition}

\theoremstyle{plain}{\theorembodyfont{\rmfamily}
}
\theoremstyle{plain}{\theorembodyfont{\rmfamily}
}
\theoremstyle{plain}{\theorembodyfont{\rmfamily}
}
\theoremstyle{plain}{\theorembodyfont{\rmfamily}
\newtheorem{example}[theorem]{Example}}

\theoremstyle{plain}{\theorembodyfont{\rmfamily}
\newtheorem{remark}[theorem]{Remark}}

\begin{document}
\begin{center}
{\bf OPTIMALITY CONDITIONS FOR A CONTROLLED SWEEPING PROCESS\\ WITH APPLICATIONS TO THE CROWD MOTION MODEL}\footnote{This research was partly supported by the National Science Foundation under grants DMS-1007132 and DMS-1512846 and by the Air Force Office of Scientific Research grant \#15RT0462.}\\[3ex]
TAN H. CAO\footnote{Department of Mathematics, Wayne State University, Detroit, Michigan 48202, USA (tan.cao@wayne.edu).} and BORIS S. MORDUKHOVICH\footnote{Department of Mathematics, Wayne State University, Detroit, Michigan 48202, USA (boris@math.wayne.edu).}
\end{center}
\small{\sc Abstract.} The paper concerns the study and applications of a new class of optimal control problems governed by a perturbed sweeping process of the hysteresis type with control functions acting in both play-and-stop operator and additive perturbations. Such control problems can be reduced to optimization of discontinuous and unbounded differential inclusions with pointwise state constraints, which are immensely challenging in control theory and prevent employing conventional variation techniques to derive necessary optimality conditions. We develop the method of discrete approximations married with appropriate generalized differential tools of modern variational analysis to overcome principal difficulties in passing to the limit from optimality conditions for finite-difference systems. This approach leads us to nondegenerate necessary conditions for local minimizers of the controlled sweeping process expressed entirely via the problem data. Besides illustrative examples, we apply the obtained results to an optimal control problem associated with of the crowd motion model of traffic flow in a corridor, which is formulated in this paper. The derived optimality conditions allow us to develop an effective procedure to solve this problem in a general setting and completely calculate optimal solutions in particular situations.
\vspace*{-0.2in}

\normalsize
\section{Introduction and Initial Discussions}
\setcounter{equation}{0}\vspace*{-0.1in}

This paper can be considered as a continuation of our work \cite{cm1}, where we formulated a new class of optimal control problems for a perturbed controlled sweeping process, justified the existence of optimal solutions therein, established the strong $W^{1,2}$-convergence of well-posed discrete approximations, and derive necessary optimality conditions for discrete optimal solutions. After summarizing the main constructions and results of \cite{cm1} in the preliminary Section~2, this paper is self-contained and may be read independently of \cite{cm1} with no additional knowledge required.\vspace*{-0.05in}

The major goal of this paper is to derive nondegenerate necessary optimality conditions for the so-called intermediate (including strong) local minimizers of the sweeping control problems under consideration by passing to the limit from the necessary optimality conditions for their discrete approximations obtained in \cite{cm1} and presented in Section~2. The class of parametric optimal control problems $(P^\tau)$ with $\tau\ge 0$ are formulated in \cite{cm1} in the following form: minimize the cost functional
\begin{equation}\label{e:4}
J[x,u,a]:=\varphi\big(x(T)\big)+\int_0^{T}\ell\big(t,x(t),u(t),a(t),\dot x(t),\dot u(t),\dot a(t)\big)dt
\end{equation}
over the control pairs $u(\cdot)\in W^{1,2}([0,T];\R^n)$ and $a(\cdot)\in W^{1,2}([0,T];\R^d)$ and the corresponding trajectories $x(\cdot)\in W^{1,2}([0,T];R^n)$ of the differential inclusion
\begin{equation}\label{e:5}
 -\dot x(t)\in N\big(x(t);C(t)\big)+f\big(x(t),a(t)\big)\;\mbox{ for a.e. }\;t\in[0,T],\quad x(0):=x_{0}\in C(0),
\end{equation}
where $x_0\in\R^n$ and $T>0$ are fixed, where the moving convex set $C(t)$ is given by
\begin{equation}\label{e:6}
C(t):=C+u(t)\;\mbox{ with }\;C:=\big\{x\in\R^n\big|\;\la x^*_i,x\ra\le 0\;\mbox{ for all }\;i=1,\ldots,m\big\}
\end{equation}
with the fixed generating $n$-vectors $x^*_i$, and where $N(x;\O)$ in \eqref{e:5} is the normal cone of convex analysis
\begin{eqnarray}\label{nor}
N(x;\O):=\big\{v\in\R^n\big|\;\la v,y-x\ra\le 0,\;y\in\O\big\}\;\mbox{ if }\;x\in\O\;\mbox{ and }\;N(x;\O):=\emp\;\mbox{ if }\;x\notin\O.
\end{eqnarray}
Besides the dynamic constraints \eqref{e:5}, problem $(P^\tau)$ involves the pointwise constraints on $u$-controls:
\begin{eqnarray}\label{e:8}
\left\{\begin{array}{ll}
\|u(t)\|=r\;\mbox{ for all }\;t\in[\tau,T-\tau],\\
r-\tau\le\|u(t)\|\le r+\tau\;\mbox{ for all }\;t\in[0,\tau)\cup(T-\tau,T]
\end{array}\right.
\end{eqnarray}
depending on the parameter $\tau\in[0,\Bar\tau]$ with $\Bar\tau:=\min\{r,T\}$ and fixed $r>0$. Note that the inclusion in \eqref{e:5} and the second part of definition \eqref{nor} implicitly yield the pointwise constraints of another type
\begin{equation}\label{mixed}
\big\la x^*_i,x(t)-u(t)\big\ra\le 0\;\mbox{ for all }\;t\in[0,T]\;\mbox{ and }\;i=1,\ldots,m.
\end{equation}\vspace*{-0.22in}

The characteristic feature of problem $(P^\tau)$ for any fixed $\tau\in[0,\Bar\tau]$ is the differential inclusion \eqref{e:5} describing, for each fixed control pair $(u(\cdot),a(\cdot))$, a perturbed version of Moreau's {\em sweeping process} \cite{mor_frict} the mathematical theory of which has been well developed; see, e.g., \cite{CT,et,KMM} and the references therein. The sweeping inclusion \eqref{e:5} significantly differs from those considered in optimal control theory for differential inclusions as developed in \cite{AS,clsw,m-book2,v} and other publications, since \eqref{e:5} admits a unique solution $x(\cdot)$ whenever the sweeping set $C(\cdot)$ and the perturbation function $a(\cdot)$ therein are given a priori; and so there is no room for optimization in such a case. Our control model in $(P^\tau)$ follows the line of \cite{chhm1,chhm2}, where control actions enter the sweeping set but not entering perturbations. Other optimal control problems for various versions of the sweeping process are considered in \cite{ao,bk,et} with no controls in the sweeping set. Namely, \cite{et} deals with controls only in perturbations addressing existence and relaxation issues for optimal solutions, while \cite{ao,bk} apply controls in associated differential equations with deriving necessary optimality conditions for discrete-time \cite{ao} and continuous-time \cite{bk} systems.\vspace*{-0.05in}

In contrast to all the previous developments, we address in this paper {\em necessary optimality conditions} for problem $(P^\tau)$ with controls in {\em both sweeping set} and {\em additive perturbations}. Note that the structure of the sweeping set in \eqref{e:5}, \eqref{e:6} is specific for the so-called {\em play-and-stop operator} \cite{smb} and largely relates to {\em rate independent hysteresis}; see, e.g., \cite{Kr,mi,smb}. Our main application here is given to a corridor version of the {\em crowd motion model} \cite{mv,ve}, where introducing controls in perturbations allows us to optimize the corresponding sweeping process and determine the optimal strategy of crowd motion participants.\vspace*{-0.05in}

Considering the triple $z=(x,u,a)\in\R^n\times\R^n\times\R^d$, it is easy to observe that \eqref{e:5} can be written as
\begin{equation}\label{e:20}
-\dot{z}(t)\in F\big(z(t)\big)\times\R^n\times\R^d\;\mbox{ a.e. }\;t\in[0,T]\;\mbox{ with }\;F(z):=N(x-u;C)+f(x,a),
\end{equation}
where the initial triple $z(0)=(x_0,u(0),a(0))$ satisfies the condition $x_0-u(0)\in C$ via the convex polyhedron $C$ defined in \eqref{e:6}. Then the sweeping optimal control problem $(P^\tau)$ amounts to minimizing the cost functional $J[z]=J[x,u,a]$ in \eqref{e:4} over $W^{1,2}$-solutions to the differential inclusion \eqref{e:20} subject to the pointwise {\em state constraints} of the equality and inequality types in \eqref{e:8} and \eqref{mixed}, where the latter ones are implicit from \eqref{e:20}. Although problem $(P^\tau)$ is now written in the usual form of the theory of differential inclusions, it is far removed from satisfying the assumptions under which necessary optimality conditions have been developed in this theory. First of all, the right-hand side of \eqref{e:20} is intrinsically {\em unbounded}, {\em discontinuous}, and {\em highly non-Lipschitzian} in any generalized sense treated by the developed approaches of optimal control theory for differential inclusions. Furthermore, besides the inequality state constraints in \eqref{mixed}, problem $(P^\tau)$ contains the unconventional {\em equality} ones as in \eqref{e:8}. Such constraints have just recently started to be considered in control theory for smooth ordinary differential equations \cite{ak}, where necessary optimality conditions are obtained under strong regularity assumptions including full rank of the smooth constraint Jacobians, which is not the case in \eqref{e:8}.\vspace*{-0.05in}

In this paper we develop the {\em method of discrete approximations} to derive necessary optimality conditions for control problems governed by differential inclusions following the scheme of \cite{m95,m-book1}, where the discrete approximation approach is realized for Lipschitzian differential inclusions without state constraints, and then its recent significant modification given in \cite{chhm2} in the case of the sweeping process with general polyhedral controlled sets but without control actions in additive perturbations. Note that our developments in this paper result in new optimality conditions that have important advantages in comparison with those in \cite{chhm2} even in the case of no controls in perturbations; see Remark~\ref{disc}(iv).\vspace*{-0.05in}

The rest of the paper is organized as follows. Section~2 summarizes some preliminary material from our preceding paper \cite{cm1} related to the strong $W^{1,2}$-convergence of discrete approximations and necessary conditions for discrete optimal solutions, which are the basic for developing here the limiting procedure to derive necessary optimality conditions in the original sweeping control problem(s) $(P^\tau)$ with the parameter $\tau$ specified above. In the main Section~3 we establish in this way, by using appropriate tools of generalized differentiation in variational analysis, necessary optimality conditions for each problem $(P^\tau)$ with $0\le\tau\le\Bar\tau$ entirely in terms of its initial data. Furthermore, we arrive at enhanced nontriviality conditions that surely exclude the appearance of the degeneracy phenomenon; cf.\ \cite{AS,v}. In Section~4 we present several numerical examples, which illustrate various specific features of the obtained optimality conditions for $(P^\tau)$ and their usage in calculating optimal solutions.\vspace*{-0.05in}

Section~5 addresses a version of the crowd motion model in a corridor that is interesting theoretically and of practical importance. We discuss this model and formulate an optimal control problem for it, which can be written in the form of the sweeping control problem studied in this paper with control functions appearing in additive perturbations. Applying necessary optimality conditions established above allows us to develop an effective procedure for calculating optimal solutions in a general setting  with finitely many participants and then fully implement it in special situations of their own interest.\vspace*{-0.05in}

Throughout the paper we use standard notation; cf.\ \cite{m-book1,rw}. Recall that $\N:=\{1,2,\ldots\}$ and that $B(x,\ve)$ stands for the closed ball of the space in question centered at $x$ with radius $\ve>0$.\vspace*{-0.15in}

\section{Discrete Approximations in Sweeping Optimal Control}
\setcounter{equation}{0}\vspace*{-0.1in}

The goal of this section is to summarize those results from our preceding paper \cite{cm1}, which make it possible to derive necessary optimality conditions for $(P^\tau)$ by passing to the limit from discrete approximations. The main assumptions made in \cite{cm1} that are standing in this paper are:

{\bf (H1)} The mapping $f\colon\R^n\times\R^d\to\R^n$ in \eqref{e:5} is continuous on $\R^n\times\R^d$ and locally Lipschitz continuous in the first argument, i.e., for every $\varepsilon>0$ there is a constant $K>0$ such that
\begin{equation*}
\left\|f(x,a)-f(y,a)\right\|\le K\left\|x-y\right\|\;\mbox{ whenever }\;(x,y)\in B(0,\varepsilon)\times B(0,\varepsilon),\quad a\in\R^d.
\end{equation*}
Furthermore, there is a constant $M>0$ ensuring the growth condition
\begin{equation}\label{gro}
\left\|f(x,a)\right\|\le M\big(1+\|x\|\big)\;\mbox{ for any }\;x\in\bigcup_{t\in[0,T]}C(t),\quad a\in\R^d.
\end{equation}\vspace*{-0.15in}

{\bf (H2)} The terminal cost function $\varphi:\R^n\to\Bar{\R}$ and the running cost function $\ell:[0,T]\times\R^{4n+2d}\to\Bar{\R}$ in \eqref{e:4} are lower semicontinuous (l.s.c.) while $\ell$ is bounded from below on bounded sets.

As follows from \cite[Theorem~1]{et} (see \cite[Proposition~2.1]{cm1}), under (H1) for any $u(\cdot)\in W^{1,2}([0,T];\R^n)$ and $a(\cdot)\in W^{1,2}([0,T];\R^d)$ there {\em exists the unique solution} $x(\cdot)\in W^{1,2}([0,T];\R^n)$ to \eqref{e:5}, \eqref{e:6} generated by $(u(\cdot),a(\cdot))$ and satisfying the {\em estimates}
\begin{equation}\label{et-es}
\left\{\begin{array}{ll}\|x(t)\|\le l:=\disp\|x_0\|+e^{2MT}\Big(2MT(1+\|x_0\|)+\int_0^T\|\dot u(s)\|ds\Big),\quad t\in[0,T],\\
\|\dot{x}(t)\|\le 2(1+l)M+\|\dot u(t)\|\;\mbox{ a.e. }\;t\in[0,T],
\end{array}\right.
\end{equation}
where the constant $M>0$ is taken from \eqref{gro}. Furthermore, it is shown in \cite[Theorem~4.1]{cm1} that if in addition to (H1) and (H2) we assume that the running cost $\ell$ in \eqref{e:4} is {\em convex} with respect to the velocity variables $(\dot x,\dot u,\dot a)$ and that $\{\dot u^k(\cdot)\}$ is bounded in $L^2([0,T];\R^n)$ while $\{a^k(\cdot)\}$ is bounded in $W^{1,2}([0,T];\R^d)$ along a minimizing sequence of $z^k(\cdot)=(x^k(\cdot),u^k(\cdot),a^k(\cdot))$ in $(P^\tau)$, then each problem $(P^\tau)$ as $\tau\in[0,\Bar\tau]$ admits an {\em optimal solution} in the class of $W^{1,2}[0,T]$ functions.\vspace*{-0.05in}

It has been realized in the conventional theory of optimal control for Lipschitzian differential inclusions (following the Bogolyubov-Young-Warga relaxation procedure in the calculus of variations and optimal control for ODEs) that the original variational problem can be replaced by its convexification with respect to velocities in such a way that the original and relaxed/convexified problems have the same values of the cost functionals and, moreover, the existing optimal solution to the relaxed problem can be approximated (in a reasonable sense) by a minimizing sequence of feasible solutions to the original one. Unfortunately, the strongest result known in this direction for non-Lipschitzian differential inclusions \cite{dfm} is not applicable to our setting in \eqref{e:20} when controls enter the sweeping set. A relaxation stability result for the perturbed sweeping process with controls acting only in perturbation is given in \cite{et}.\vspace*{-0.05in}

Mentioning this, denote by $\ell_{F}(t,x,u,a,\dot x,\dot u,\dot a)$ the {\em convexification} of the integrand in \eqref{e:4} on the set $F(x,u,a)$ from \eqref{e:20} with respect to $(\dot x,\dot u,\dot a)$, i.e., the largest convex and l.s.c.\ function majorized by $\ell(t,x,u,a,\cdot,\cdot,\cdot)$ on this set for all $t,x,u,a$ with putting $\Hat\ell:=\infty$ at points out of the set $F(x,u,a)$, and define the {\em variational relaxation} of $(P^\tau)$ as follows:
\begin{equation*}
\mbox{minimize }\;\Hat{J}[z]:=\varphi\big(x(T)\big)+\int^{T}_{0}\Hat{\ell}_{F}\big(t,x(t),u(t),a(t),\dot{x}(t),\dot{u}(t),\dot{a}(t)\big)dt
\end{equation*}
over all $z(\cdot)=(x(\cdot),u(\cdot),a(\cdot))\in W^{1,2}[0,T]$ satisfying \eqref{e:8}. Similarly to \cite{m95}, we say that $\bar{z}(\cdot)$ is a {\em relaxed intermediate local minimizer} (r.i.l.m) for $(P^\tau)$ if it is feasible to this problem with $J[\bar{z}]=\Hat{J}[\bar{z}]$ and if there are numbers $\alpha\ge 0$ and $\epsilon>0$ such that $J[\bar{z}]\le J[z]$ for any feasible solution $z(\cdot)$ to $(P^\tau)$ with
\begin{equation}\label{ilm}
\left\|z(t)-\oz(t)\right\|<\epsilon\;\mbox{ for all }\;t\in[0,T]\;\mbox{ and }\;\alpha\int^{T}_{0}\left\|\dot{z}(t)-\dot{\oz}(t)\right\|^{2}dt<\epsilon.
\end{equation}
The word ``relaxed" can be skipped here if the the running cost $\ell$ is convex in $(\dot x,\dot u,\dot a)$, or--more generally--under the local relaxation stability of $(P^\tau)$. The intermediate local minimum defined is obviously situated between the classical notions of weak and strong minima. Furthermore, we do not restrict the generality by putting $\al=1$ in \eqref{ilm} while deriving necessary optimality conditions.

The main goal of this paper is to establish necessary optimality conditions for the {\em given r.i.l.m.}\ $\oz(\cdot)=(\ox(\cdot),\ou(\cdot),\oa(\cdot))$ of problem $(P^\tau)$. To accomplish it, we use the method of discrete approximations. Proceeding as in \cite{cm1}, for all $k\in\N$ consider the discrete mesh
\begin{equation}\label{mesh}
\Delta_{k}:=\big\{0=t^{k}_{0}<t^{k}_{1}<\ldots<t^{k}_{k}\big\}\;\mbox{ with }\;h_{k}:=t^{k}_{j+1}-t^{k}_{j}\dn 0\;\mbox{ as }\;k\to\infty
\end{equation}
on $[0,T]$, denote by $j_\tau(k):=\left[k\tau/T\right]$ the smallest index $j$ with $t^k_j\ge\tau$ and by $j^\tau(k):=[k(T-\tau)/T]-1$ the largest $j$ with $t^k_j\le T-\tau$. Then we construct the discrete-time optimization problems $(P^\tau_k)$ as follows. Given $\ve>0$ from \eqref{ilm} together with the numbers $\Tilde\mu>0$ and $\ve_k\dn 0$ as $k\to\infty$ explicitly defined in \cite[Theorem~3.1]{cm1} via $\oz(\cdot)$ and the data of $(P^\tau)$, each discrete problem $(P^\tau_k)$ consists of minimizing the cost
\begin{eqnarray*}
\begin{aligned}
J_{k}[z^{k}]:=&\varphi(x^{k}_{k})+h_{k}\sum^{k-1}_{j=0}\ell\left(t^{k}_{j},x^{k}_{j},u^{k}_{j},a^{k}_{j},\dfrac{x^{k}_{j+1}-x^{k}_{j}}{h_{k}},\dfrac{u^{k}_{j+1}
-u^{k}_{j}}{h_{k}},\dfrac{a^{k}_{j+1}-a^{k}_{j}}{h_{k}}\right)\\
+&\sum^{k-1}_{j=0}\int^{t^{k}_{j+1}}_{t^{k}_{j}}{\left(\left\|\dfrac{x^{k}_{j+1}-x^{k}_{j}}{h_{k}}-\dot{\bar{x}}(t)\right\|^{2}+\left\|\dfrac{u^{k}_{j+1}
-u^{k}_{j}}{h_{k}}-\dot{\bar{u}}(t)\right\|^{2}+\left\|\dfrac{a^{k}_{j+1}-a^{k}_{j}}{h_{k}}-\dot{\bar{a}}(t)\right\|^{2}\right)}dt\\
+&\dist^2\left(\left\|\dfrac{u^k_1-u^k_0}{h_k}\right\|;\big(-\infty,\Tilde{\mu}\big]\right)+\dist^2\left(\sum^{k-2}_{j=0}\left\|\dfrac{u^k_{j+2}-2u^k_{j+1}
+u^k_j}{h_k}\right\|;\big(-\infty,\Tilde{\mu}\big]\right)
\end{aligned}
\end{eqnarray*}
over elements $z^{k}:=(x^{k}_{0},x^{k}_{1},\ldots,x^{k}_{k},u^{k}_{0},u^{k}_{1},\ldots,u^{k}_{k-1},a^{k}_{0},a^{k}_{1},\ldots,a^{k}_{k-1})$ satisfying the constraints
\begin{equation}\label{disc-dyn}
x^{k}_{j+1}\in x^{k}_{j}-h_{k}F(x^{k}_{j},u^{k}_{j},a^{k}_{j})\;\mbox{ for }\;j=0,\ldots,k-1\;\mbox{ with }\;(x^{k}_{0},u^k_0,a^k_0)=\big(x_{0},\ou(0),\oa(0)\big),
\end{equation}
\begin{equation*}
\left\langle x^*_{i},x^{k}_k-u^{k}_k\right\rangle\le 0\;\mbox{ for }\;i=1,\ldots,m,
\end{equation*}
\begin{equation*}
\|u^k_j\|=r\;\mbox{for}\;j=j_\tau(k),\ldots,j^\tau(k);\;r-\tau-\varepsilon_k\le\|u^k_j\|\le r+\tau+\varepsilon_k\;\mbox{for}\;j\le j_\tau(k)-1\;\mbox{and}\;j\ge j^\tau(k)+1,
\end{equation*}
\begin{equation*}
\left\|(x^{k}_{j},u^{k}_{j},a^{k}_{j})-\big(\bar{x}(t^{k}_{j}),\bar{u}(t^{k}_{j}),\bar{a}(t^{k}_{j})\big)\right\|\le\epsilon/2\;\mbox{ for }\;j=0,\ldots,k-1,
\end{equation*}
\begin{equation*}
\sum^{k-1}_{j=0}\int^{t^{k}_{j+1}}_{t^{k}_{j}}{\left(\left\|\dfrac{x^{k}_{j+1}-x^{k}_{j}}{h_{k}}-\dot{\bar{x}}(t)\right\|^{2}+\left\|\dfrac{u^{k}_{j+1}-u^{k}_{j}}
{h_{k}}-\dot{\bar{u}}(t)\right\|^{2}+\left\|\dfrac{a^{k}_{j+1}-a^{k}_{j}}{h_{k}}-\dot{\bar{a}}(t)\right\|^{2}\right)}dt\le\dfrac{\epsilon}{2},
\end{equation*}
\begin{equation*}
\left\|\dfrac{u^k_1-u^k_0}{t^k_1-t^k_0}\right\|\le\Tilde{\mu}+1,\;\mbox{ and }\;\sum^{k-2}_{j=0}\left\|\dfrac{u^k_{j+2}-2u^k_{j+1}+u^k_j}{h_k}\right\|\le \Tilde{\mu}+1,
\end{equation*}\vspace*{-0.05in}

Denote further the collection of {\em active constraint indices} of polyhedron \eqref{e:6} at $\ox\in C$ by
\begin{equation}\label{e:18}
I(\ox):=\big\{i\in\{1,\ldots,m\}\big|\;\big\la x^*_i,\ox\big\ra=0\big\},
\end{equation}
and recall the explicit representation of the set $F(z)$ in \eqref{e:20} that follows from \cite[Proposition~3.1]{hmn}:
\begin{equation}\label{e:19}
F(z)=\Big\{\sum_{i\in I(x-u)}\lambda_{i}x^*_i\Big|\;\lambda_i\ge 0\Big\}+f(x,a)\;\mbox{ with }\;x-u\in C.
\end{equation}
We also need to recall the {\em featured} subsets of the active constraint indices \eqref{e:18} defined by
\begin{equation}\label{feat}
I_0(y):=\big\{i\in I(\ox)\big|\;\left\langle x^*_i,y\right\rangle=0\big\}\;\mbox{ and }\;I_>(y):=\big\{i\in I(\ox)\big|\;\left\langle x^*_i,y\right\rangle >0\big\},\quad y\in\R^n.
\end{equation}\vspace*{-0.22in}

Now we are ready to summarize the main results of \cite{cm1} (cf.\ Theorems~4.1, 5.2, and 7.2) in one theorem used in what follows. Note that although $(P^\tau_k)$ is {\em intrinsically nonsmooth} (primarily because of the discrete dynamics \eqref{disc-dyn}, even in the case of smooth cost functions $\ph$ and $\ell$ that is not assumed here), the necessary conditions for discrete optimal solutions presented below do not contain any generalized normals and derivatives associated with $F$. This is due to the {\em explicit calculation} of the employed construction in terms of the initial sweeping data, which is given in \cite[Theorem~6.1]{cm1} on the basic of \cite{hmn} and calculus rules. The only construction of generalized differentiation needed in the following theorem and also employed below is the {\em subdifferential} of a function $\psi\colon\R^n\to\oR:=(-\infty,\infty]$ finite at $\os$ that is defined by
\begin{equation}\label{sub}
\partial\psi(\os):=\Big\{v\in\R^n\Big|\;\exists\,s_k\st{\psi}{\to}\os,\;v_k\to v\;\mbox{ with }\;\liminf_{k\to\infty}\frac{\psi(s_k)-\psi(\os)-\la v_k,s_k-s\ra}{\|s_k-\os\|}\ge 0\Big\},
\end{equation}
where the symbol $s_k\st{\psi}{\to}\os$ indicates that $s_k\to\os$ and $\psi(s_k)\to\psi(\os)$. It reduces to the classical derivative for smooth functions to the subdifferential of convex analysis if $\psi$ is convex. For the general class of l.s.c.\ functions, the subdifferential \eqref{sub} enjoys comprehensive calculus rules based of variational/extremal principles of variational analysis; see \cite{m-book1,rw} for more details and references. The local smoothness of the perturbation function $f$ is imposed below for simplicity.\vspace*{-0.05in}

In the following theorem we impose the assumptions of \cite[Theorem~3.1]{cm1} on the given r.i.l.m.\ $\oz(\cdot)$ for $(P^\t)$ satisfying \eqref{e:20} at all $t_j$ as $j=0,\ldots,k-1$ (with the right-side derivative at $t_0=0$) while do not formulate them for brevity. Recall only that these assumptions hold if $\oz(\cdot)\in W^{2,\infty}[0,T]$.\vspace*{-0.1in}

\begin{theorem}{\bf(strong convergence and necessary conditions for discrete optimal solutions).}\label{Th:11}
Let $\bar{z}(\cdot)=(\bar{x}(\cdot),\bar{u}(\cdot),\bar{a}(\cdot))$ be a r.i.l.m.\ for problem $(P^\tau)$ with $0\le\tau\le\Bar\tau:=\min\{r,T\}$. In addition to $(H1)$, $(H2)$ and all the assumptions of {\rm\cite[Theorem~3.1]{cm1}} imposed on $\oz(\cdot)$, suppose that the cost functions $\ph$ and $\ell(t,\cdot,\cdot)$ are locally Lipschitzian around $\ox(T)$ and $(\oz(t),\dot\oz(t))$ for a.e.\ $t\in[0,T]$, respectively, that $\ell(\cdot,z,\dot z)$ is a.e.\ continuous on $[0,T]$ being uniformly majorized by a summable function near $\oz(\cdot)$, and that $f(\cdot,\cdot)$ from \eqref{e:20} is smooth around $(\ox(t),\oa(t))$ on $[0,T]$. Then each problem $(P^\tau_k)$ for $k\in\N$ sufficiently large admits an optimal solution $\oz^k(\cdot)=(\ox^k(\cdot),\ou^k(\cdot),\oa^k(\cdot))$ such that their piecewise linear extensions on $[0,T]$ converge to $\oz(\cdot)$ as $k\to\infty$ in the norm topology of $W^{1,2}[0,T]$ with
\begin{equation}\label{e:40}
\disp\lim_{k\to\infty}\frac{\ox^k_1-x^k_0}{h_k}=\dot\ox(0),\quad\limsup_{k\to\infty}\left\|\dfrac{\ou^k_1-\ou^k_0}{h_k}\right\|\le\Tilde{\mu},\quad\limsup_{k\to\infty}\sum^{k-2}_{j=0}
\left\|\dfrac{\ou^k_{j+2}-2\ou^k_{j+1}+\ou^k_j}{h_k}\right\|\le\Tilde{\mu}.
\end{equation}
Define further the vectors $(\theta^{xk}_j,\theta^{uk}_j,\theta^{ak}_j)\in\R^n\times\R^n\times\R^d$ for $j=0,\ldots,k-1$ by
\begin{equation}\label{e:61}
(\theta^{xk}_j,\theta^{uk}_j,\theta^{ak}_j):=2\int^{t_{j+1}}_{t_j}\bigg(\dfrac{\ox^{k}_{j+1}-\ox^{k}_{j}}{h_k}-\dot{\ox}(t),\dfrac{\ou^{k}_{j+1}-\ou^{k}_{j}}{h_k}
-\dot{\ou}(t),\dfrac{\oa^{k}_{j+1}-\oa^{k}_{j}}{h_k}-\dot{\oa}(t)\bigg)dt,
\end{equation}
the vectors $\chi^k_j\in\R^n$ for $j=0,\ldots,k-1$ by
\begin{equation}\label{chi}
\chi^k_j:=\left\{\begin{array}{ll}
\disp\frac{x^k_1-x^k_0}{h_k}&\mbox{if }\;j=1,\\
0&\mbox{otherwise},
\end{array}\right.
\end{equation}
and assume that the generating elements $\{x^*_i|\;i\in I(\ox(t)-\ou(t)\}$ of \eqref{e:6} along the corresponding active indices \eqref{e:18} are linearly independent for all $t\in[0,T]$. Then there exist triples $(\lm^k,\xi^k,p^k)$ with $\lm_k\ge 0$, $\xi^{k}=(\xi^{k}_{0},\ldots,\xi^{k}_{k})\in\R^{k+1}$, and $p^k_j=(p^{xk}_j,p^{uk}_j,p^{ak}_j)\in\R^n\times\R^n\times\R^d$ as $j=0,\ldots,k$ together with vectors $\eta^k_j\in\R^m_+$ as $j=0,\ldots,k$ and $\gamma^k_{j}\in\R^m$ as $j=0,\ldots,k-1$, and subgradients
\begin{equation}\label{e:60}
\big(w^{xk}_{j},w^{uk}_{j},w^{ak}_{j},v^{xk}_{j},v^{uk}_{j},v^{ak}_{j}\big)\in\partial\ell\bigg(\oz^{k}_{j},\dfrac{\oz^{k}_{j+1}-\oz^{k}_{j}}{h_k}\bigg)\;\mbox{ for }\;j=0,\ldots,k-1,
\end{equation}
where we drop indicating the time-dependence of the running cost $\ell$ in \eqref{e:4} for brevity, such that the following necessary optimality conditions are satisfied:\vspace*{-0.05in}

$\bullet$ {\sc Nontriviality condition}
\begin{equation}\label{e:73}
\lm^k+\left\|\xi^k\right\|+\left\|p^{uk}_0\right\|+\left\|p^{ak}_0\right\|\not=0.
\end{equation}\vspace*{-0.15in}

$\bullet$ {\sc Primal-dual dynamic relationships}
\begin{equation}\label{e:63}
\dfrac{\ox^k_{j+1}-\ox^k_j}{-h_k}-f(\ox^k_j,\oa^k_j)=\sum_{i\in I(\ox^k_j-\ou^k_j)}\eta^k_{ji}x^*_i,
\end{equation}
\begin{equation}\label{e:66}
\begin{array}{ll}
\dfrac{p^{xk}_{j+1}-p^{xk}_j}{h_k}-\lm^kw^{xk}_j-\chi^k_j=\nabla_xf(\ox^k_j,\oa^k_j)^*\left(\lm^k(v^{xk}_j+h^{-1}_k\theta^{xk}_j)-p^{xk}_{j+1}\right)\\
+\disp\sum_{i\in I_0\left(-p^{xk}_{j+1}+\lm^k(h^{-1}_k\theta^k_{xj}+v^{xk}_j)\right)\cup I_>\left(-p^{xk}_{j+1}+\lm^k(h^{-1}_k\theta^k_{xj}+v^{xk}_j)\right)}\gamma^k_{ji}x^*_i,
\end{array}
\end{equation}
\begin{equation}\label{e:67}
\begin{array}{ll}
\dfrac{p^{uk}_{j+1}-p^{uk}_j}{h_k}-\lm^kw^{uk}_j-\dfrac{2}{h_k}\xi^k_j\ou^k_j=-\disp\sum_{i\in I_0\left(-p^{xk}_{j+1}+\lm^k(h^{-1}_k\theta^k_{xj}+v^{xk}_j)\right)\cup I_>\left(-p^{xk}_{j+1}+\lm^k(h^{-1}_k\theta^k_{xj}+v^{xk}_j)\right)}\gamma^k_{ji}x^*_i,
\end{array}
\end{equation}
\begin{equation}\label{e:68}
\dfrac{p^{ak}_{j+1}-p^{ak}_j}{h_k}-\lm^kw^{ak}_j=\nabla_a f(\ox^k_j,\oa^k_j)^*\left(\lm^k(v^{xk}_j+h^{-1}_k\theta^{xk}_j)-p^{xk}_{j+1}\right),
\end{equation}
\begin{equation}\label{e:69}
\big[\left\la x^*_i,\ox^k_j-\ou^k_j\right\ra<0\big]\Longrightarrow\eta^k_{ji}=0,\quad\eta^k_{ji}>0\Longrightarrow\big[\la x^*_i,-p^{xk}_{j+1}+\lm^k\left(h^{-1}_k\theta^k_{xj}+v^{xk}_j\right)\ra=0\big],
\end{equation}
\begin{equation}\label{gamma-k}
\big[\la x^*_i,\ox^k_j-\ou^k_j\ra<0\big]\Longrightarrow\gamma^k_{ji}
\end{equation}
valid for all the indices $j=0,\ldots,k-1$ and $i=1,\ldots,m$, where the featured index subsets $I_0(\cdot)$ and $I_>(\cdot)$ are taken from \eqref{feat}. Furthermore, we have
\begin{equation}\label{e:65}
\left\{\begin{array}{lcl}
i\in I_0\left(-p^{xk}_{j+1}+\lm^k(h^{-1}_k\theta^k_{xj}+v^{xk}_j)\right)\Longrightarrow\gamma^k_{ji}\in\R,\\
i\in I_>\left(-p^{xk}_{j+1}+\lm^k(h^{-1}_k\theta^k_{xj}+v^{xk}_j)\right)\Longrightarrow\gamma^k_{ji}\ge 0,\\
\big[i\not\in I_0\left(-p^{xk}_{j+1}+\lm^k(h^{-1}_k\theta^k_{xj}+v^{xk}_j)\right)\cup I_>\left(-p^{xk}_{j+1}+\lm^k(h^{-1}_k\theta^k_{xj}+v^{xk}_j)\right)\big] \Longrightarrow\gamma^k_{ji}=0,
\end{array}\right.
\end{equation}
\begin{equation}\label{e:70}
\xi^k_j\in N\big(\|\ou^k_j\|;[r-\t-\ve_k,r+\t+\ve_k]\big)\;\mbox{ for }\;j\in\big\{0,\ldots,j_r(k)-1\big\}\cup\big\{j^k(k)+1,\ldots,k\big\}.
\end{equation}\vspace*{-0.15in}

$\bullet$ {\sc Transversality conditions} at the right endpoint:
\begin{equation}\label{e:72}
-p^{xk}_k\in\lm^k\partial\varphi(\ox^k_k)+\sum^m_{i=1}\eta^k_{ki}x^*_i,\quad p^{uk}_k=\sum^m_{i=1}\eta^k_{ki}x^*_i-2\xi^k_k\ou^k_k,\quad p^{ak}_k=0
\end{equation}
together with the endpoint implications
\begin{equation}\label{e:71}
\big[\la x^*_i,\ox^k_k-\ou^k_k\ra<0\big]\Longrightarrow\eta^k_{ki}=0\;\mbox{ for }\;i=1,\ldots,m.
\end{equation}
\end{theorem}\vspace*{-0.2in}

\section{Necessary Conditions for Sweeping Optimal Solutions}
\setcounter{equation}{0}\vspace*{-0.1in}

In this section we derive necessary optimality conditions for relaxed intermediate local minimizers of the sweeping control problem $(P^\tau)$ under consideration in the general case of $0\le\tau\le\Bar\t=\min\{r,T\}$ with some specifications and improvements in the case where $\t$ is not an endpoint. For convenience in the formulation of some conditions in the next theorem, the symbol $N(x;\O)$ is used therein for the normal cone to the convex set in question at $x\in\O$ instead of its explicit description given in \eqref{nor}. For the same reason, the symbol $D^*G$ is used to indicate the coderivative of the sweeping set-valued mapping appeared in the theorem despite its explicit calculation in terms of the initial sweeping data presented below.\vspace*{-0.05in}

To specify this issue, recall that the {\em coderivative} of an arbitrary set-valued mapping $G\colon\R^n\tto\R^m$ at the point of the graph $(\ox,\oy)\in\gph G:=\{(x,y)\in\R^n\times\R^m|\;y\in G(x)\}$ can be defined by
\begin{equation}\label{cod}
D^*G(\ox,\oy)(u):=\big\{v\in\R^n\big|\;(v,-u)\in\partial\dd\big((\ox,\oy);\gph G\big)\big\},\quad u\in\R^m,
\end{equation}
via the subdifferential \eqref{sub} of the indicator function $\dd(\cdot;\gph G)$ of the graphical set $\gph G$ that equals to 0 on $\gph G$ and to $\infty$ otherwise. Recall also that $\dom G:=\{x\in\R^n|\;G(x)\ne\emp\}$. If $G$ is single-valued and smooth around $\ox$, this construction reduces to $\{\nabla G(\ox)^*u\}$ via the adjoint/transpose Jacobian matrix $\nabla G(\ox)$ of $G$ at $\ox$, while in general it is set-valued, nonconvex, and enjoys full calculus; see \cite{m-book1}.\vspace*{-0.05in}

The set-valued mapping $G\colon\R^n\tto\R^n$ of our main interest here, which appears in Theorem~\ref{Th:12} and other results below, is the {\em normal cone mapping} $G(x):=N(x;C)$ generated by the underlying convex polyhedron $C$ from \eqref{e:6}. The coderivative of this mapping is known as the {\em second-order subdifferential} of $\dd(\cdot;C)$ in the sense of \cite{m-book1} while being precisely calculated in terms of the initial data of $C$ in the following proposition, which is taken from \cite[Theorem~4.6]{hmn}.\vspace*{-0.1in}

\begin{proposition} {\bf (precise calculating the coderivatives of the normal cone mappings associated with convex polyhedra).}\label{Th:7} Let $G(x)=N(x;C)$ be the normal cone mapping associated with the convex polyhedron \eqref{e:6}, and let the featured active index subsets $I_0(\cdot)$ and $I_>(\cdot)$ be defined in \eqref{feat}. Given $(\ox,\oy)\in\gph G$, assume that the generating elements $\{x^*_i|\;i\in I(\bar{x})\}$ of \eqref{e:6} along the active constraint indices \eqref{e:18} are linearly independent. Then we have the coderivative expression
\begin{equation*}
D^*G(\ox,\oy)(u)=\span\big\{x^*_i\big|\;i\in I_0(u)\big\}+\cone\big\{x^*_i\big|\;i\in I_>(u)\big\}\;\mbox{ for all }\;u\in\dom D^*G(\ox,\oy),
\end{equation*}
where the latter coderivative domain is characterized by
\begin{equation*}
u\in\dom D^*G(\ox,\oy)\Longleftrightarrow\big[i\in J(\ox,\oy)\Longrightarrow\la x^*_i,u\ra=0\big]
\end{equation*}
via the so-called strict complementarity subset of active indices $J(\ox,\oy):=\{i\in I(\ox)|\;\lm_i>0\}$ for a unique collection of the multipliers $\lm_i\ge 0$ coming from the representation $\oy=\sum_{i\in I(\ox)}\lm_i x^*_i$.
\end{proposition}\vspace{-0.1in}

Before establishing the aforementioned necessary optimality conditions for $(P^\tau)$ we make the following remark. It is well known that the subdifferential mapping \eqref{sub} used in condition \eqref{e:60} of Theorem~\ref{Th:11} is {\em robust} (i.e., closed-graph) with respect to subdifferentiation variables. However, in the discrete-time and continuous-time settings under consideration allows the dependence of the running cost $\ell$ on the time parameter, which is not under subdifferentiation. For the limiting procedure in what follows, we {\em require} the subdifferential robustness with respect to the time parameter. It is not a restrictive assumption that holds, in particular, for smooth functions with time-continuous derivatives as well as in rather general nonsmooth settings discussed, e.g., in \cite{m95,m-book2}. We recall also that (as assumed in \cite[Theorem~3.1]{cm1}) that the local optimal solution $\oz(\cdot)$ under consideration satisfies the sweeping differential inclusion \eqref{e:20} at all the the mesh points with the right and left derivatives at $t=0$ and $t=T$, respectively.\vspace*{-0.1in}

\begin{theorem}{\bf (general necessary optimality conditions for the perturbed sweeping process).}\label{Th:12} Let $\oz(\cdot)=(\ox(\cdot),\ou(\cdot),\oa(\cdot))$ be a r.i.l.m.\ for problem $(P^\tau)$ with any $\t\in[0,\Bar\t]$. In addition to the assumptions of Theorem~{\rm\ref{Th:11}}, suppose that $\ell$ in \eqref{e:4} is continuous in $t$ a.e.\ on $[0,T]$ and admits the representation
\begin{equation}\label{separ}
\ell(t,z,\dot z)=\ell_1(t,z,\dot x)+\ell_2(t,\dot u)+\ell_3(t,\dot a),
\end{equation}
where the $($local$)$ Lipschitz constants of $\ell_1(t,\cdot,\cdot)$ and $\ell_3(t,\cdot)$ can be chosen as essentially bounded on $[0,T]$ and continuous at a.e.\ $t\in[0,T]$ including $t=0$, while $\ell_2$ is differentiable in $\dot u$ on $\R^n$ satisfying
\begin{equation}\label{estimate1}
\begin{array}{c}
\|\nabla_{\dot{u}}\ell_2(t,\dot{u},\dot{a})\|\le L\|\dot{u}\|\;\mbox{ and }\;\|\nabla_{\dot{u}}\ell_2(t,\dot{u}_1)-\nabla_{\dot{u}}\ell_2(t,\dot{u}_2)\|\le L|t-s|+L\|\dot{u}_1-\dot{u}_2\|
\end{array}
\end{equation}
for some constant $L>0$, all numbers $t,s\in[0,T]$, $\dot a\in\R^d$, and all vectors $\dot u,\dot u_1,\dot u_2\in\R^n$.\vspace*{-0.05in}

Then there exist a number $\lm\ge 0$, an adjoint arc $p(\cdot)=(p^x(\cdot),p^u(\cdot),p^a(\cdot))\in W^{1,2}([0,T];\R^n\times\R^n\times\R^d)$, subgradient functions $w(\cdot)=(w^x(\cdot),w^u(\cdot),w^a(\cdot))\in L^2([0,T];\R^{2n+d})$ and $v(\cdot)=(v^x(\cdot,v^u(\cdot),v^a(\cdot))\in L^2([0,T];\R^{2n+d})$ well defined at $t=0$  and satisfying the inclusion
\begin{equation}\label{subg}
\big(w(t),v(t)\big)\in\co\partial\ell\big(t,\oz(t),\dot{\oz}(t)\big)\;\mbox{ for a.e. }\;t\in[0,T],
\end{equation}
and measures $\gamma=(\gg_1,\ldots,\gg_n)\in C^*([0,T];\R^n)$, $\xi\in C^*([0,T];\R)$ on $[0,T]$ such that we have:

$\bullet$ {\sc Primal-dual dynamic relationships:}
\begin{equation}\label{e:76}
-\dot{\ox}(t)=\sum^m_{i=1}\eta_i(t)x^*_i+f\big(\ox(t),\oa(t)\big)\;\mbox{ for a.e. }\;t\in[0,T],
\end{equation}
where $\eta_i(\cdot)\in L^2([0,T];\R_+)$ are well defined at $t=T$ being uniquely determined by representation \eqref{e:76};
\begin{equation}\label{e:77}
\dot{p}(t)=\lm w(t)+\Big(\nabla_x f\big(\ox(t),\oa(t)\big)^*\big(\lm v^x(t)-q^x(t)\big),0,\nabla_a f\big(\ox(t),\oa(t)\big)^*\big(\lm v^x(t)-q^x(t)\big)\Big),
\end{equation}
\begin{equation}\label{e:78}
q^u(t)=\lm\nabla_{\dot u}\ell\big(t,\dot{\ou}(t)\big),\quad q^a(t)\in\lm\partial_{\dot a}\ell_3\big(t,\dot{\oa}(t)\big)\;\mbox{ for a.e. }\;t\in[0,T],
\end{equation}
where the vector function $q=(q^x,q^u,q^a):[0,T]\to\R^n\times\R^n\times\R^d$ is of bounded variation with its left-continuous representative given for all $t\in[0,T]$, except at most a countable subset, by
\begin{equation}\label{e:79}
q(t):=p(t)-\disp\int_{[t,T]}\big(d\gamma(s),2\ou(s)d\xi(s)-d\gamma(s),0\big).
\end{equation}
Moreover, for a.e.\ $t\in[0,T]$ including $t=T$ and all $i=1,\ldots,m$ we have the implications
\begin{equation}\label{e:74}
\la x^*_i,\ox(t)-\ou(t)\ra<0\Longrightarrow\eta_i(t)=0,\quad\eta_i(t)>0\Longrightarrow\la x^*_i,\lm v^x(t)-q^x(t)\ra=0.
\end{equation}\vspace*{-0.22in}

$\bullet$ {\sc Transversality conditions} at the right and left endpoints, respectively:
\begin{equation}\label{e:81}
\left\{\begin{array}{ll}
-p^x(T)-\disp\sum_{i\in I(\ox(T)-\ou(T))}\eta_i(T)x^*_i\in\lm\partial\ph\big(\ox(T)\big),\\
p^u(T)-\disp\sum_{i\in I(\ox(T)-\ou(T))}\eta_i(T)x^*_i\in 2\ou(T)N\big(\|\ou(T)\|;[r-\t,r+\t]\big),\quad p^a(T)=0;
\end{array}\right.
\end{equation}
\begin{equation}\label{e:81a}
\left\{\begin{array}{ll}
q^x(0)\in\R^n,\quad q^a(0)=\lm v^a(0),\quad q^u(0)\in\lm v^u(0)-2\ou(0)N\big(\|\ou(0)\|;[r-\t,r+\t]\big)\\\\
+D^*G\big(x_0-\ou(0),-\dot{\ox}(0)-f(\ox(0),\oa(0)\big)\big(-q^x(0)+\lm v^x(0)\big)
\end{array}\right.
\end{equation}
with the coderivative $D^*G$ of $G(\cdot)=N(\cdot;C)$ explicitly calculated in Proposition~{\rm\ref{Th:7}} and with
\begin{equation}\label{e:81b}
\disp\sum_{i\in I(\ox(T)-\ou(T))}\eta_i(T)x^*_i\in N\big(\ox(T)-\ou(T);C\big).
\end{equation}\vspace*{-0.2in}

$\bullet$ {\sc Measure nonatomicity conditions:}\\[1ex]
{\bf(a)} If $t\in[0,T)$ and $\la x^*_i,\ox(t)-\ou(t)\ra<0$ for all $i=1,\ldots,m$, then there exists a neighborhood $V_t$ of $t$ in $[0,T)$ such that $\gamma(V)=0$ for all Borel subsets $V$ of $V_t$.\\[1ex]
{\bf(b)} Assume that $\t\in(0,\Bar\t)$ and take any $t\in[0,\tau)\cup(T-\tau,T]$ with $r-\tau<\|\ou(t)\|< r+\tau$. Then there exists a neighborhood $W_t$ of $t$ in $(0,\tau)\cup(T-\tau,T)$ such that $\xi(W)=0$ for all Borel subsets $W$ of $W_t$.

$\bullet$ {\sc Nontriviality conditions:}\\[1ex]
{\bf(a)} Impose one of the following assumptions on the local minimizer $\oz(\cdot)$ and the data of $(P^\tau)$ as $\t\in[0,\Bar\t]$:
\begin{equation}\label{assum-nontr}
\mbox{either }\;l(r+2\t)<(r-2\t)^2,\quad\mbox{or}\quad\la\ox(t),\ou(t)\ra\not=\|\ou(t)\|^2\;\mbox{ for all }\;t\in[0,T),
\end{equation}
where the constant $l>0$ is calculated in \eqref{et-es} with $u(\cdot)=\ou(\cdot)$.\footnote{Note that the first condition in \eqref{assum-nontr} implies the second one for $\t=0$, while in general they are independent.} Then we have
\begin{equation}\label{e:83}
\lm+\|q^u(0)\|+\|p(T)\|>0
\end{equation}
provided that either $\t<r$ or $\ou(T)\ne 0$.\\[1ex]
{\bf (b)} If in addition $0<\t<r$, the we have the following enhanced nontriviality conditions while imposing the corresponding endpoint interiority assumptions:
\begin{equation}\label{enh1}
\big[\la x^*_i,x_0-\ou(0)\ra<0,\;r-\tau<\|\ou(0)\|<r+\tau,\;i=1,\ldots,m\big]\Longrightarrow\big[\lm+\|p(T)\|>0\big].
\end{equation}
\begin{equation}\label{enh2}
\big[\la x^*_i,\ox(T)-\ou(T)\ra<0,\;r-\tau<\|\ou(T)\|<r+\tau,\;i=1,\ldots,m\big]\Longrightarrow\big[\lm+\|q^u(0)\|>0\big].
\end{equation}
\end{theorem}\vspace*{-0.1in}
{\bf Proof.} The derivation of the necessary optimality conditions for the given r.i.l.m.\ $\oz(\cdot)$ in problem $(P^\t)$ is based on passing the limit as $k\to\infty$ from the optimality conditions for the strongly convergent sequence $\oz^k(\cdot)\to\oz(\cdot)$ of optimal solutions to the discrete problems $(P^\t_k)$ obtained in Theorem~\ref{Th:11}. The proof is rather involved, and for the reader's convenience we split it into several steps.\\[1ex]
{\bf Step~1:} {\em Subdifferential inclusion.} Let us first justify \eqref{subg}. For each $k\in\N$ define the functions $w^k,v^k\colon[0,T]\to\R^{2n+d}$ as piecewise constant extensions to $[0,T]$ of the vectors $w^k_j$ and $v^k_j$ that are defined on the mesh $\Delta_k$ and satisfy the subdifferential inclusion \eqref{e:60} therein. The assumptions made and the structure of $\ell$ in \eqref{separ}, \eqref{estimate1} ensure that the subgradient sets $\partial\ell(t,\cdot)$ are uniformly bounded near $\oz(\cdot)$ by the $L^2$-Lipschitz constant of $\ell$, and thus the sequence $\{(w^k(\cdot),v^k(\cdot))\}$ is weakly compact in $L^2([0,T];\R^{2(2n+d)}):=L^2[0,T]$. This allows us to select a subsequence (no relabeling hereafter) converging
$$
\big(w^k(\cdot),v^k(\cdot)\big)\to\big(w(t),v(t)\big)\;\mbox{ weakly in }\;L^2[0,T]\;\mbox{ as }\;k\to\infty
$$
to some $(w(\cdot),v(\cdot))\in L^2[0,T]$. Furthermore, the local Lipschitz continuity of $\ell(0,\cdot,\cdot)$ yields by \eqref{e:60} for $j=0$ that the sequence $\{(w^k_0,v ^k_0)\}$ is bounded and hence converges as $k\to\infty$ to a pair $(w_0,v_0)=:(w(0),v(0))$ along a subsequence. It follows from the aforementioned Mazur weak closure theorem that there are convex combinations of $(w^k(\cdot),v^k(\cdot))$, which converge to $(w(\cdot),v(\cdot))$ in the $L^2$-norm and hence a.e.\ on $[0,T]$ for some subsequence. Then passing to the limit in \eqref{e:60} along the latter subsequence and taking into account the assumed a.e.\ continuity of the running cost $\ell$ in $t$ and robustness of its subdifferential in $(z,\dot z)$ with respect to all the variables, we arrive at the convexified inclusion \eqref{subg}.\\[1ex]
{\bf Step~2:} {\em Passing to the limit in the primal equation.} Our next aim is to arrive at the primal equation \eqref{e:76} and the first implication in \eqref{e:74} with the corresponding functions $\eta_i(\cdot)$ by passing to the limit in \eqref{e:63} and \eqref{e:69}. We start with considering the functions
$$
\theta^{k}(t):=\dfrac{\theta^{k}_j}{h_k}\;\mbox{ as }\;t\in[t^k_j,t^k_{j+1}),\;j=0,\ldots,k-1,\quad k\in\N,
$$
on $[0,T]$ with $\theta^{k}_j$ taken from \eqref{e:61}. It follows from the convergence $\oz^k(\cdot)\to\oz(\cdot)$ in Theorem~\ref{Th:11} that
\begin{equation}\label{theta}
\begin{aligned}
\int^T_0\|\theta^{xk}(t)\|^2dt=&\sum^{k-1}_{j=0}\dfrac{\|\theta^{xk}_j\|^2}{h_k}\le\dfrac{4}{h_k}\sum^{k-1}_{j=0}\bigg(\int^{t^k_{j+1}}_{t^k_j}\left\|\dot{\ox}(t)
-\dfrac{\ox^k_{j+1}-\ox^k_j}{h_k}\right\|dt\bigg)^2\\\le&4\sum^{k-1}_{j=0}\bigg(\int^{t^k_{j+1}}_{t^k_j}\left\|\dot{\ox}(t)-\dfrac{\ox^k_{j+1}
-\ox^k_j}{h_k}\right\|^2dt=4\int^T_0\left\|\dot{\ox}(t)-\dot{\ox}^k(t)\right\|^2dt\to 0
\end{aligned}
\end{equation}
and similarly for $\theta^{uk}(\cdot)$ and $\theta^{ak}(\cdot)$. This yields the a.e.\ convergence of these functions to zero on $[0,T]$. Moreover, the construction above shows that we can always have $\theta^k_0\to 0=:\th(0)$. \vspace*{-0.05in}

Further, it is easy to see that the assumed linear independence of $\{x^*_i\left|\right.i\in I(\ox(\cdot)-\ou(\cdot))\}$ ensures the one for $\{x^*_i\left|\right.i\in I(\ox^k_j-\ou^k_j)\}$ by definition \eqref{e:18} and the strong convergence of Theorem~\ref{Th:11}. This allows us to take the vectors $\eta^k_j\in\R^m_+$ from Theorem~\ref{Th:11} and construct the piecewise constant functions $\eta^k(\cdot)$ on $[0,T]$ by $\eta^k(t):=\eta^k_j$ for $t\in[t^k_j,t^k_{j+1})$ with $\eta^k(T):=\eta^k_k$. It follows from \eqref{e:65} that
\begin{equation}\label{e:85}
-\dot{\ox}^k(t)=\sum^m_{i=1}\eta^k_i(t)x^*_i+f\big(\ox^k(t^k_j),\oa^k(t^k_j)\big)\;\mbox{ whenever }\;t\in(t^k_j,t^k_{j+1}),\quad k\in\N,
\end{equation}
via the corresponding components of $\eta^k(t)$. On the other hand, the feasibility of $\oz(\cdot)$ to $(P^\tau)$ yields $-\dot{\ox}(t)\in G(\ox(t)-\ou(t))+f(\ox(t),\oa(t))$ for a.e. $t\in[0,T]$ with the closed-valued normal cone mapping $G(\cdot)=N(\cdot;C)$. Due to the measurability of $G(\cdot)$ by \cite[Theorem~4.26]{rw} and the measurable selection result from \cite[Corollary~4.6]{rw} we find nonnegative measurable functions $\eta_i(\cdot)$ on $[0,T]$ as $i=1,\ldots,m$ such that equation \eqref{e:76} and the first implication in \eqref{e:74} hold. Combining \eqref{e:85} and \eqref{e:76} gives us
$$
\dot{\ox}(t)-\dot{\ox}^k(t)=\sum^m_{i=1}\big[\eta^k_i(t)-\eta_i(t)\big]x^*_i+f\big(\ox^k(t^k_j),\oa^k(t^k_j)\big)-f\big(\ox(t),\oa(t)\big)
$$
for $t\in(t^k_j,t^k_{j+1})$ and $j=0,\ldots,k-1$. Thus we get the estimate
\begin{equation}\label{e:86}
\left\|\sum^m_{i=1}\big[\eta_i(t)-\eta^k_i(t)\big]x^*_i\right\|\le\|\dot{\ox}(t)-\dot{\ox}^k(t)\|+\big\|f\big(\ox(t),\oa(t)\big)-f\big(\ox^k(t^k_j),\oa^k(t^k_j)
\big)\big\|
\end{equation}
for $t\in(t^k_j,t^k_{j+1})$. For each $k\in\N$ define now the function
\begin{equation}\label{e:87}
\nu^k(t):=\max\big\{t^k_j\big|\;t^k_j\le t,\;0\le j\le k\big\}\;\mbox{ for all }\;t\in[0,T].
\end{equation}
Passing to the limit in \eqref{e:86} with replacing $t^k_j$ by $\nu(t)$ and taking into account the strong convergence $\oz^k(\cdot)\to\oz(\cdot)$ together with the continuity of $f$ on the left-hand side of \eqref{e:86}, we get
$$
\sum_{i\in I(\ox(t)-\ou(t))}\big[\eta_i(t)-\eta^k_i(t)\big]x^*_i\to 0\;\mbox{ as }\;k\to\infty\;\mbox{ for a.e. }\;t\in[0,T].
$$
Then the assumed linear independence of the generating vectors $x^*_i$ with $i\in I(\ox(t)-\ou(t))$ ensures the a.e.\ convergence $\eta^k(t)\to\eta(t)$ on $[0,T]$ as $k\to\infty$. Furthermore, we can always get that $\eta^k_k$ converge to the well-defined vector $(\eta_1(T),\ldots,\eta_m(T))$. Proceeding similarly to the proof of \cite[Theorem~6.1]{chhm2}, we can justify the extra regularity $\eta(\cdot)\in L^2([0,T];\R_+^m)$, which however is not used in what follows.\\[1ex]
{\bf Step~3:} {\em Extensions of approximating dual elements.} Here we extend discrete dual elements from Theorem~\ref{Th:11} to the whole interval $[0,T]$. First construct $q^k(t)=(q^{xk}(t),q^{uk}(t),q^{ak}(t))$ on $[0,T]$ as the piecewise linear extensions of $q^k(t^k_j):=p^k_j$ as $j=0,\ldots,k$. Then define $\gamma^k(t)$ on $[0,T]$ as piecewise constant $\gamma^k(t):=\gamma^k_j$ for $t\in[t^k_j,t^k_{j+1})$ and $j=0,\ldots,k-1$ with $\gamma^k(t^k_k):=0$. We also set $\xi^k(t):=\dfrac{\xi^k_j}{h_k}$ for $t\in[t^k_j,t^k_{j+1})$ and $j=0,\ldots,k-1$ with $\xi^k(t^k_k):=\xi^k_k$. Using the function $\nu^k(t)$ given in \eqref{e:87}, we deduce respectively from \eqref{e:66}, \eqref{e:67}, and \eqref{e:68} that
\begin{equation}\label{e:88}
\begin{aligned}
\dot{q}^{xk}(t)-\lm^kw^{xk}(t)=&\nabla_xf\big(\ox^k(\nu^k(t)),\oa^k(\nu^k(t)\big)^*\big(\lm^k(v^{xk}(t)+\theta^{xk}(t))-q^{xk}(\nu^k(t)+h_k)\big)\\
+&\sum_{i\in I_0(\lm^k(v^{xk}(t)+\theta^{xk}(t))-q^{xk}(\nu^k(t)+h_k))\cup I>(\lm^k(v^{xk}(t)+\theta^{xk}(t))-q^{xk}(\nu^k(t)+h_k))}\gamma^k_i(t)x^*_i,
\end{aligned}
\end{equation}
\begin{equation}\label{e:89}
\begin{aligned}
\dot{q}^{uk}(t)-\lm^kw^{uk}(t)=&2\xi^k(t)\ou^k\big(\nu^k(t)\big)\\
-&\sum_{i\in I_0(\lm^k(v^{xk}(t)+\theta^{xk}(t))-q^{xk}(\nu^k(t)+h_k))\cup I>(\lm^k(v^{xk}(t)+\theta^{xk}(t))-q^{xk}(\nu^k(t)+h_k))}\gamma^k_i(t)x^*_i,
\end{aligned}
\end{equation}
\begin{equation}\label{e:90}
\dot{q}^{ak}(t)-\lm^kw^{ak}(t)=\nabla_af\big(\ox^k(\nu^k(t)),\oa^k(\nu^k(t)\big)^*\big(\lm^k(v^{xk}(t)+\theta^{xk}(t))-q^{xk}(\nu^k(t)+h_k)\big)
\end{equation}
for $t\in(t^k_j,t^k_{j+1})$ and $j=0,\ldots,k-1$. Next we define $p^k(t)=(p^{xk}(t),p^{uk}(t),p^{ak}(t))$ on $[0,T]$ by setting
\begin{equation}\label{e:91}
p^k(t):=q^k(t)+\int_{[t,T]}\left(\sum^m_{i=1}\gamma_i^k(s)x^*_i,2\xi^k(s)\ou^k\big(\nu^k(s)\big)-\sum^m_{i=1}\gamma^k_i(s)x^*_i,0\right)ds
\end{equation}
for all $t\in[0,T]$. This gives us $p^k(T)=q^k(T)$ and the differential relation
\begin{equation}\label{e:92}
\dot{p}^k(t)=\dot{q}^k(t)-\left(\sum^m_{i=1}\gamma^k_i(t)x^*_i,2\xi^k(t)\ou^k\big(\nu^k(t)\big)-\sum^m_{i=1}\gamma^k_i(t)x^*_i,0\right)\;\mbox{ a.e. }\;t\in[0,T].
\end{equation}
It follows from \eqref{e:92}, \eqref{e:88}--\eqref{e:90}, and the definition of $I_0(\cdot)$ and $I_>(\cdot)$ in \eqref{feat} that
\begin{equation}\label{e:93}
\dot{p}^{xk}(t)-\lm^kw^{xk}(t)-\chi^k_j=\nabla_xf\big(\ox^k(\nu^k(t)),\oa^k(\nu^k(t)\big)^*\big(\lm^k(v^{xk}(t)+\theta^{xk}(t))-q^{xk}(\nu^k(t)+h_k)\big),
\end{equation}
\begin{equation}\label{e:94}
\dot{p}^{uk}(t)-\lm^kw^{uk}(t)=0,
\end{equation}
\begin{equation}\label{e:95}
\dot{p}^{ak}(t)-\lm^kw^{ak}(t)=\nabla_a f\big(\ox^k(\nu^k(t)),\oa^k(\nu^k(t)\big)^*\big(\lm^k(v^{xk}(t)+\theta^{xk}(t))-q^{xk}(\nu^k(t)+h_k)\big),
\end{equation}
for every $t\in(t^k_j,t^k_{j+1})$,\;$j=0,\ldots,k-1$. Define now the vector measures $\gamma^k_{mes}$ and $\xi^k_{mes}$ on $[0,T]$ by
\begin{equation}\label{e:96}
\int_A d\gamma^k_{mes}:=\int_A\sum^m_{i=1}\gamma^k_i(t)x^*_idt,\quad\int_Ad\xi^k_{mes}:=\int_A\xi^k(t)dt\;\mbox{ for any Borel subset }\;A\subset[0,T]
\end{equation}
with dropping further the symbol ``$mes$" for simplicity. By taking into account the preservation of all the relationships in Theorem~\ref{Th:11} by normalization and the above constructions of the extended functions on $[0,T]$, we can rewrite the nontriviality condition \eqref{e:73} as
\begin{equation}\label{e:97}
\lm^k+\|p^k(T)\|+\|q^{uk}(0)\|+\|q^{ak}(0)\|+\int^T_0|\xi^k(t)|dt+|\xi^k_k|+\int^T_0\left\|\sum^m_{i=1}\gamma^k_i(t)x^*_i\right\|dt=1,\quad k\in\N.
\end{equation}
{\bf Step~4:} {\em Passing to the limit in dual dynamic relationships.} Using \eqref{e:97} allows us to suppose without loss of generality that $\lm^k\to\lm$  as $k\to\infty$ for some $\lm\ge 0$. Let us next verify that the sequence $\{(p^{xk}_0,\ldots,p^{xk}_k)\}_{k\in\N}$ is bounded in $\R^{(k+1)n}$. Indeed, we have by \eqref{e:66} that
\begin{equation}\label{pxk}
p^{xk}_j=p^{xk}_{j+1}-\lm^kh_kw^{xk}_j-h_k\chi^k_j-\nabla_xf(\ox^k_j,\oa^k_j)^*(\lm^kh_kv^{xk}_j+\lm^k\theta^{xk}_j-h_kp^{xk}_{j+1})-h_k\sum^m_{i=1}\gamma^k_{ji}x^*_i
\end{equation}
for $j=0,\ldots,k-1$. It follows from \eqref{theta} and \eqref{e:97} that the quantities $\nabla_xf(\ox^k_j,\oa^k_j)$, $\lm^k\theta^{xk}_j$, and $h_k\sum^m_{i=1}\gamma^k_{ji}x^*_i$ are uniformly bounded for $j=0,\ldots,k-1$ while $\chi^k_j\to 0$ as $k\to\infty$ due to definition \eqref{chi} and the first condition in \eqref{e:40}. Furthermore, the imposed structure \eqref{separ} of $\ell$ and the assumptions on the Lipschitz constant $L(t)$ of the running cost in \eqref{e:4}, which are equivalent to the Riemann integrability of $L(\cdot)$ on $[0,T]$, yield by \eqref{subg} the relationships
\begin{equation}\label{tL}
\|h_kw^{xk}_j\|=\|h_kw^{xk}(t_j)\|\le h_k L(t_j)\le\sum^{k-1}_{i=0}L(t_i)h_k\le 2\int_{[0,T]}L(t)dt=:\Tilde L<\infty
\end{equation}
and ensure similarly that $\|h_kv^{xk}_j\|<\Tilde L$. Thus we get from \eqref{pxk} that the boundedness of $\{p^{xk}_j\}_{k\in\N}$ follows from the boundedness of $\{p^{xk}_{j+1}\}_{k\in\N}$. Since the sequence of $p^k_k=p^k(T)$ is bounded by \eqref{e:97}, we therefore justify the claim on the boundedness of $\{(p^{xk}_0,\ldots,p^{xk}_k)\}_{k\in\N}$.\vspace*{-0.05in}

To deal with the functions $q^{xk}(\cdot)$, we derive from their construction and the equations in \eqref{e:66} that
\begin{equation}\label{qk}
\begin{aligned}
\sum^{k-1}_{j=0}\|q^{xk}(t_{j+1})-q^{xk}(t_j)\|&\le\lm^k\sum^{k-1}_{j=0}h_kw^{xk}_j+h_k\sum^{k-1}_{j=0}\big\|\nabla_xf\big(\ox^k_j,\oa^k_j\big)^*\big(\lm^k(\theta^{xk}(t_j)
-p^{xk}_{j+1})\big)\big\|\\
&+\sum^{k-1}_{j=0}\|\nabla_x f(\ox^k_j,\oa^k_j)^*(\lm^kh_kv^{xk}_j)\|+\int^T_0\left\|\sum^m_{i=1}\gamma^k_i(t)x^*_i\right\|dt.
\end{aligned}
\end{equation}
It comes from \eqref{tL} that the first term on the right-hand side of \eqref{qk} is bounded by $\lm_k\Tilde L$. We also have
$$
h_k\sum^{k-1}_{j=0}\|\nabla_xf\big(\ox^k_j,\oa^k_j\big)^*\big(\lm^k(\theta^{xk}(t_j)-p^{xk}_{j+1})\big)\|\le T\max_{0\le j\le k-1} \bigg\{\big\|\nabla_xf\big(\ox^k_j,\oa^k_j\big)^*\big(\lm^k(\theta^{xk}(t_j)-p^{xk}_{j+1})\big)\|\bigg\},
$$
which ensures the boundedness of the second term on the right-hand side of \eqref{qk} by the boundedness of $\{p^{xk}_j\}_{k\in\N}$. Similarly we get the boundedness of the third term on the right-hand side of \eqref{qk}, while this property of the forth term therein follows from \eqref{e:97}. This shows by estimate \eqref{qk} and the construction of $q^{xk}(t)$ on $[0,T]$ that the functions $q^{xk}(\cdot)$ are of uniformly bounded variation on this interval. In the same way we verify that $q^{uk}(\cdot)$ and $q^{ak}(\cdot)$ are of uniformly bounded variation on $[0,T]$ and arrive therefore at this conclusion for the whole triple $q^k(\cdot)$. It clearly implies that
\begin{eqnarray*}
2\|q^{k}(t)\|-\|q^{k}(0)\|-\|q^{k}(T)\|\le\|q^{k}(t)-q^{k}(0)\|+\|q^{k}(T)-q^{k}(t)\|\le\var(q^{k};[0,T])\;\mbox{ for all }\;t\in[0,T],
\end{eqnarray*}
which justifies the uniform boundedness of $q^k(\cdot)$ on $[0,T]$ is since both sequences $\{q^{k}(0)\}$ and $\{q^{k}(T)\}$ are bounded by \eqref{e:97}. Then the classical Helly selection theorem allows us to find a function of bounded variation $q(\cdot)$ such that $q^k(t)\to q(t)$ as $k\to\infty$ pointwise on $[0,T]$. Employing further \eqref{e:97} and the measure construction in \eqref{e:96} tell us that the measure sequences $\{\gamma^k\}$ and $\{\xi^k\}$ are bounded in $C^*([0,T];\R^n)$ and $C^*([0,T];\R)$ respectively. It follows from the weak$^*$ sequential compactness of the unit balls in these spaces that there are measures $\gamma\in C^*([0,T];\R^n)$ and $\xi\in C^*([0,T];\R)$ such that the pair $(\gamma^k,\xi^k)$ weak$^*$ converges to $(\gamma,\xi)$ along some subsequence.\vspace*{-0.05in}

Combining the uniform boundedness of $q^k(\cdot)$, $w^k(\cdot)$, and $v^k(\cdot)$ on $[0,T]$ with \eqref{e:91}, \eqref{e:93}--\eqref{e:95}, and \eqref{e:97} allows us to deduce that the sequence $\{p^k(\cdot)\}$ is bounded in $W^{1,2}([0,T];\R^{3n})$ and hence weakly compact in this space. By the Mazur weak closure theorem we find $p(\cdot)\in W^{1,2}([0,T];\R^n)$ such that a sequence of convex combinations of $\dot{p}^k(t)$ converges to $\dot{p}(t)$ for a.e.\ $t\in[0,T]$. Passing now to the limit in \eqref{e:93}--\eqref{e:95} as $k\to\infty$ and using \eqref{theta}, we arrive at the representation of $\dot p(\cdot)$ in \eqref{e:77}.
\vspace*{-0.05in}

Next we proceed with deriving adjoint relationships involving the limiting function $q(\cdot)$ of bounded variation on $[0,T]$. Note to this end that if $\eta_i(t)>0$ for some $t\in[0,T]$ and $i\in\{1,\ldots,m\}$, then $\eta_i^k(t)>0$ whenever $k$ is sufficiently large due to the a.e.\ convergence $\eta_i^k(\cdot)\to\eta_i(\cdot)$ on $[0,T]$. This implies by \eqref{e:69} that $\la x^*_i,-q^{xk}(\nu(t)+h_k)+\lm^k(\theta^{xk}(t)+v^{xk}(t))\ra=0$ for such $k$ and $t$, and so we arrive at $\la x^*_i,-q^x(t)+\lm v^x(t)\ra=0$ while $k\to\infty$, which thus justifies the second implication in \eqref{e:74}.\vspace*{-0.05in}

Remembering the construction of $q^k(\cdot)$ in Step~3 allows us to rewrite \eqref{e:67} and \eqref{e:68} as, respectively,
\begin{equation}\label{e:98}
q^{uk}\big(\nu(t)+h_k\big)=\lm^k\big(v^{uk}(t)+\theta^{uk}(t)\big)\;\mbox{ and }\;q^{ak}\big(\nu(t)+h_k\big)=\lm^k\big(v^{ak}(t)+\theta^{ak}(t)\big)
\end{equation}
for $t\in(t^k_j,t^k_{j+1})$ and $j=0,\ldots,k-1$. Passing to the limit in \eqref{e:98} with taking into account \eqref{subg} and the assumptions on $\ell_2,\ell_3$ in \eqref{separ}, we arrive at both equations in \eqref{e:78}. Observe further that
$$
\left\|\int_{[t,T]}\sum^m_{i=1}\gamma^k_i(s)x^*_ids-\int_{[t,T]}d\gamma(s)\right\|=\left\|\int_{[t,T]}d\gamma^k(s)-\int_{[t,T]}d\gamma(s)\right\|\to 0\;\mbox{ as}\;k\to\infty
$$
for all $t\in[0,T]$ except a countable subset of $[0,T]$ by the weak$^*$ convergence of the measures $\gamma^k$ to $\gamma$ in the space $C^*([0,T];\R^n)$; cf.\ \cite[p.\ 325]{v} for similar arguments. This ensures by \eqref{e:96} that
\begin{equation}\label{e:99}
\int_{[t,T]}\sum^m_{i=1}\gamma^k_i(s)x^*_ids\to \int_{[t,T]}d\gamma(s)\;\mbox{ as }\;k\to\infty.
\end{equation}
To obtain \eqref{e:79} by passing to the limit in \eqref{e:92}, consider next the estimate
\begin{equation}\label{e:100}
\begin{aligned}
&\left\|\int_{[t,T]}\xi^k(s)\ou^k\big(\nu^k(s)\big)ds-\int_{[t,T]}\ou(s)d\xi(s)\right\|\\
\le&\left\|\int_{[t,T]}\xi^k(s)\ou^k\big(\nu^k(s)\big)ds-\int_{[t,T]}\xi^k(s)\ou(s)ds\right\|
+\left\|\int_{[t,T]}\xi^k(s)\ou(s)ds-\int_{[t,T]}\ou(s)d\xi(s)\right\|\\
=&\left\|\int_{[t,T]}\xi^k(s)\big[\ou^k\big(\nu^k(s)\big)-\ou(s)\big]ds \right\|+\left\|\int_{[t,T]}\xi^k(s)\ou(s)ds-\int_{[t,T]}\ou(s)d\xi(s)\right\|
\end{aligned}
\end{equation}
and observe that the first summand in the rightmost part of \eqref{e:100} disappears as $k\to\infty$  due to the uniform convergence $\ou^k(\cdot)\to \ou(\cdot)$ on $[0,T]$ and the uniform boundedness of $\int^T_0|\xi^k(t)|dt$ by \eqref{e:97}. The second summand therein also converges to zero for all $t\in[0,T]$ except some countable subset by the weak$^*$ convergence $\xi^k\to\xi$ in $C^*([0,T];\R)$. Hence we get
$$
\int_{[t,T]}\xi^k(s)\ou^k(\tau^k(s))ds\to\int_{[t,T]}\ou(s)d\xi(s)\;\mbox{ as }\;k\to\infty
$$
and thus arrive at \eqref{e:79} by passing to the limit in \eqref{e:92}. Finally at this step, observe that the implications is \eqref{e:74} follow directly by passing to the limit in their discrete counterparts \eqref{e:69} and \eqref{e:71}.\\[1ex]
{\bf Step~5:} {\em Transversality conditions.} Let us first verify the right endpoint one \eqref{e:81}. For all $k\in\N$ we have by the second condition in \eqref{e:72} and the normal cone representation from \eqref{e:19} that
\begin{equation}\label{tr1}
p^{uk}_k+2\xi^k_k\ou^k_k=\sum^m_{i=1}\eta^k_{ki}x^*_i=\sum_{i\in I(\ox^k_k-\ou^k_k)}\eta^k_{ki}x^*_i\in N(\ox^k_k-\ou^k_k;C),
\end{equation}
where $\eta^k_{ki}=0$ for $i\in\{1,\ldots,m\}\backslash I(\ox^k_k-\ou^k_k)$. Denote $\vt_k:=\sum_{i\in I(\ox^k_k-\ou^k_k)}\eta^k_{ki}x^*_i$ and observe that a subsequence of $\{\vt_k\}$ converges to some $\vt$ due to the boundedness of $\{\xi^k_k\}$ by \eqref{e:97} and the convergence of $\{p^{uk}_k\}$ and $\{\ou^k_k\}$ established above. Furthermore, it follows from the robustness of the normal cone in \eqref{tr1}, the convergence $\ox^k_k-\ou^k_k\to\ox(T)-\ou(T)$, and the inclusion $I(\ox^k_k-\ou^k_k)\subset I(\ox(T)-\ou(T))$ for large $k\in\N$ that $\vt\in N(\ox(T)-\ou(T);C)$. Similarly the inclusion in \eqref{e:72} tells us that
\begin{equation}\label{tr2}
-p^{xk}_k-\vt_k\in\lm^k\partial\ph(\ox^k_k)\;\mbox{ for all }\;k\in\N.
\end{equation}
Passing now to the limit as $k\to\infty$ in \eqref{tr1}, \eqref{tr2}, inclusion \eqref{e:70} for $\xi^k_k$, and the second condition in \eqref{e:72} with taking into account the robustness of the subdifferential in \eqref{tr2}, we arrive at the relationships
\begin{equation*}
-p^x(T)-\vt\in\lm\partial\varphi\big(\ox(T)\big),\quad p^u(T)-\vt\in-2\ou(T)N\big(\|\ou(T)\|;[r-\tau,r+\tau]\big),\quad p^a(T)=0
\end{equation*}
with $\vt=\sum_{i\in I(\ox(T)-\ou(T))}\eta_i(T)x^*_i\in N(\ox(T)-\ou(T);C)$. This clearly verifies the transversality conditions at the right endpoint in \eqref{e:81} supplemented by \eqref{e:81b}.\vspace*{-0.05in}

To justify the left endpoint transversality \eqref{e:81a}, we deduce from \eqref{e:67} and \eqref{e:68} for $j=0$ as well as the conditions on $\gg^k_0$ in
Theorem~\ref{Th:11} and the coderivative definition \eqref{cod} that
$$
p^{uk}_0+h_k\lm^kw^{uk}_0+2\xi^k_0\ou^k_0-\lm^k(v^{uk}_0+\theta^{uk}_0)\in D^*G\Big(\ox^k_0-\ou^k_0,-\dfrac{\ox^k_1-\ox^k_0}{h_k}-f(\ox^k_0,\oa^k_0)\Big)\big(-p^{xk}_1+\lm^k(\theta^{xk}_0+v^{xk}_0)\big),
$$
$$p^{ak}_0=\lm^k(v^{ak}_0+\theta^{ak}_0-h_k\lm^k w^{ak}_0-h_k\nabla_a f(\ox^k_0,\oa^k_0)^*\big(\lm^k(\theta^{xk}_0+v^{xk}_0)\;\mbox{ whenever }\;k\in\N.
$$
Now we can pass to the limit as $k\to\infty$ in these relationships by taking into account \eqref{e:70} for $j=0$, the first condition in \eqref{e:40}, the construction of $q^k(t^k_j)=p^k_j$, the convergence statements for $w^k_0$, $v^k_0$, $\th^k_0$ established above as well as robustness of the normal cone and coderivative. This readily gives us \eqref{e:81a}.\\[1ex]
{\bf Step~6:} {\em Measure nonatomicity conditions.} To verify condition (a) therein without any restriction on $\t\in[0,\Bar\t]$, pick $t\in[0,T)$ with $\la x^*_i,\ox(t)-\ou(t)\ra<0$ for $i=1,\ldots,m$ and by continuity of $(\ox(\cdot),\ou(\cdot))$ find a neighborhood $V_t$ of $t$ such that $\la x^*_i,\ox(s)-\ou(s)\ra<0$ whenever $s\in V_t$ and $i=1,\ldots,m$. This shows by the established convergence of the discrete optimal solutions that $\la x^*_i,\ox^k(t^k_j)-\ou^k(t^k_j)\ra<0$ if $t^k_j\in V_t$ for $i=1,\ldots,m$ for all $k\in\N$ sufficiently large. Then it follows from \eqref{e:69} that $\gamma^k(t)=0$ on any Borel subset $V$ of $V_t$, and therefore $\|\gamma^k\|(V)=\int_Vd\|\gamma^k\|=\int_V\|\gamma^k(t)\|dt=0$ by the construction of the measure $\gg^k$ in \eqref{e:96}. Passing now to limit as $k\to\infty$ and taking into account the measure convergence obtained in Step~3, we arrive at $\|\gamma\|(V)=0$ and thus justify the first measure nonatomicity condition. The proof of the nonatomicity condition (b) for the measure $\xi$ is similar provided the choice of $\t\in(0,\Bar\t)$.\\[1ex]
{\bf Step~7:} {\em General nontriviality condition.} Let us justify the nontriviality condition \eqref{e:83} for any $\t\in[0,\Bar\t]$ under the assumptions made therein. Suppose on the contrary to \eqref{e:83} that $\lm=0$, $q^u(0)=0$, and $p(T)=0$, which yields $\lm^k\to 0$, $q^{uk}(0)\to 0$, and $p^{k}(T)\to 0$ as $k\to\infty$. It follows from \eqref{e:79} that $q^a(\cdot)=p^a(\cdot)$ is absolutely continuous on $[0,T]$, which allows us to deduce from \eqref{e:78} and the established convergence of adjoint trajectories that $p^{ak}_0=q^{ak}(0)\to 0$ as $k\to\infty$. Furthermore, we have by the construction of $\xi^k(\cdot)$ on $[0,T]$ in Step~3 the equalities
$$
\int^T_0|\xi^k(t)|dt=\sum^{k-1}_{j=0}h_k\dfrac{|\xi^k_j|}{h_k}=\sum^{k-1}_{j=0}|\xi^k_j|.
$$
Taking into account that $I_0\left(-p^{xk}_{j+1}+\lm^k(h^{-1}_k\theta^k_{xj}+v^{xk}_j)\right)\cup I_>\left(-p^{xk}_{j+1}+\lm^k(h^{-1}_k\theta^k_{xj}+v^{xk}_j)\right)\subset I(\ox^k_j-\ou^k_j)$ and that $\la x^*_i,\ox^k_j-\ou^k_j\ra =0$ for all $i\in I(\ox^k_j-\ou^k_j)$ and $j=0,\ldots,k-1$, it follows from \eqref{e:67} that
\begin{equation}\label{e:102}
2\xi^k_j\la\ou^k_j,\ox^k_j-\ou^k_j\ra=\la p^{uk}_{j+1}-p^{uk}_j-h_k\lm^kw^{uk}_j,\ox^k_j-\ou^k_j\ra,\quad j=0,\ldots,k-1.
\end{equation}
Using now the first condition in \eqref{assum-nontr} imposed on the initial data of $(P^\t)$ and that $\|\ou^k_j\|=r$ for all $j=j_\tau(k),\ldots,k$ while $r-\tau-\varepsilon_k\le\|\ou^k_j\|\le r+\tau+\varepsilon_k$ for $j=0,\ldots,j_\tau(k)-1$ and small $\varepsilon_k<\tau$, we get
$$
|\la\ou^k_j,\ox^k_j\ra|\le\|\ou^k_j\|\cdot\|\ox^k_j\|\le(r+2\tau)l<(r-2\tau)^2<r^2=\|\ou^k_j\|^2,
$$
which immediately implies the validity of the estimate
\begin{equation}\label{u-est}
|\la\ou^k_j,\ox^k_j-\ou^k_j\ra|=\big|\la\ou^k_j,\ox^k_j\ra-\|\ou^k_j\|^2\big|>0
\end{equation}
whenever $k$ is sufficiently large. On the other hand, \eqref{u-est} follows directly from the alternative assumption in \eqref{assum-nontr} imposed on the optimal solution to $(P^\t)$. Employing \eqref{u-est} and the equalities in \eqref{e:100} gives us
\begin{equation}\label{e:103}
2|\xi^k_j|\le\left(\left\|p^{uk}_{j+1}-p^{uk}_j\right\|+h_k\lm^k\|w^{uk}_j\|\right)\dfrac{\left\|\ox^k_j-\ou^k_j\right\|}{|\la\ou^k_j,\ox^k_j-\ou^k_j\ra|},\quad j=0,\ldots,k-1.
\end{equation}
The obvious boundedness of $\left\{\dfrac{\left\|\ox^k_j-\ou^k_j\right\|}{|\la\ou^k_j,\ox^k_j-\ou^k_j\ra|}\right\}$ allows us to assume without loss of generality that $\dfrac{\left\|\ox^k_j-\ou^k_j\right\|}{|\la\ou^k_j,\ox^k_j-\ou^k_j\ra|}\le 1$ for $j=0,\ldots,k-1$, and then we get from \eqref{e:103} that
\begin{equation}\label{e:104}
2\sum^{k-1}_{j=0}|\xi^k_j|\le \sum^{k-1}_{j=0}\|p^{uk}_{j+1}-p^{uk}_j\|+h_k\lm^k\sum^{k}_{j=0}\|w^{uk}_j\|
\end{equation}
The second sum in \eqref{e:104} disappears as $k\to\infty$ due to the assumptions on $\ell_1$; see \eqref{tL} in Step~4. To proceed with the first sum in \eqref{e:104}, we have the estimates
\begin{equation}\label{e:105}
\begin{aligned}
\sum^{k-1}_{j=0}\|p^{uk}_{j+1}-p^{uk}_j\|\le&\sum^{k-1}_{j=1}\|p^{uk}_{j+1}-p^{uk}_j\|+\|p^{uk}_1\|+\|p^{uk}_0\|\\
\le&\lm^k\sum^{k-1}_{j=1}\left\|\dfrac{\theta^{uk}_j-\theta^{uk}_{j-1}}{h_k}\right\|+\lm^k\dfrac{\|\theta^{uk}_0\|}{h_k}+\sum^{k-1}_{j=1}\|v^{uk}_j-v^{uk}_{j-1}\|
+\lm^k\|v^{uk}_0\|+\|p^{uk}_0\|\\\le&2\lm^k\sum^{k-1}_{j=1}\left\|\dfrac{\ou^k_{j+1}-2\ou^k_j+\ou^k_{j-1}}{h_k}\right\|+2\lm^k\sum^{k-1}_{j=1}\left\|
\dfrac{\ou(t^k_{j+1})-2\ou(t^k_j)+\ou(t^k_{j-1})}{h_k}\right\|\\
+&\lm^k\dfrac{\|\theta^{uk}_0\|}{h_k}+\sum^{k-1}_{j=1}\|v^{uk}_j-v^{uk}_{j-1}\|+\lm^k\|v^{uk}_0\|+\|p^{uk}_0\|\\
\le&4\Tilde\mu\lm^k+\lm^k \dfrac{\|\theta^{uk}_0\|}{h_k}+\sum^{k-1}_{j=1}\|v^{uk}_j-v^{uk}_{j-1}\|+\lm^k\|v^{uk}_0\|+\|p^{uk}_0\|,
\end{aligned}
\end{equation}
where $\Tilde\mu$ is defined in Theorem~\ref{Th:11}. The running cost structure \eqref{separ} and differentiability of $\ell_2$ in $\dot{u}$ yield
$$
v^{uk}_j=\nabla_{\dot{u}}\ell_2\bigg(t_j,\dfrac{\ou^k_{j+1}-\ou^k_j}{h_k},\dfrac{\oa^k_{j+1}-\oa^k_j}{h_k}\bigg)\;\mbox{ for }\;j=0,\ldots,k-1.
$$
Then the third estimate in \eqref{estimate1} ensures that
\begin{equation*}
\sum^{k-1}_{j=1}\|v^{uk}_j-v^{uk}_{j-1}\|\le\sum^{k-1}_{j=1}L\bigg(\big(t_{j+1}-t_j\big)+\bigg\|\dfrac{\ou^k_{j+1}-2\ou^k_j+\ou^k_{j-1}}{h_k}\bigg\|\bigg)\le L(T+\Tilde{\mu}).
\end{equation*}
Deduce further from the definition of $\th^{uk}$ in \eqref{e:61} the representation
$$
\dfrac{\theta^{uk}_0}{h_k}=\dfrac{2(\ou^k_1-\ou^k_0)}{h_k}-\dfrac{2\big(\ou(h_k)-\ou(0)\big)}{h_k}
$$
and observe that $\lm^k\dfrac{\|\theta^{uk}_0\|}{h_k}\to 0$ as $k\to\infty$ due to second estimate in \eqref{e:40} and the assumption imposed on $\ou(\cdot)$ in Theorem~\ref{Th:11} via \cite[Theorem~3.1]{cm1}. We have furthermore that
$$
\|v^{uk}_0\|=\bigg\|\nabla_{\dot{u}}\ell_2\bigg(0,\dfrac{\ou^k_1-\ou^k_0}{h_k},\dfrac{\oa^k_1-\oa^k_0}{h_k}\bigg)\bigg\|\le L\bigg\|\dfrac{\ou^k_1-\ou^k_0}{h_k}\bigg\|\le L\Tilde{\mu}
$$
due to \eqref{e:40} and the second estimate in \eqref{estimate1}. This shows therefore that
\begin{equation}\label{e:105a}
\sum^{k-1}_{j=0}\|p^{uk}_{j+1}-p^{uk}_j\|\to 0\;\mbox{ and }\int^T_0|\xi^k(t)|dt\to 0\;\mbox{ as }\;k\to\infty
\end{equation}
by \eqref{e:105} and \eqref{e:104}, respectively. Appealing again to \eqref{e:67} gives us
$$\int^T_0\left\|\sum^m_{i=1}\gamma^k_i(t)x^*_i\right\|dt=\sum^{k-1}_{j=0}\left\|h_k\sum^m_{i=1}\gamma^k_{ji}x^*_i\right\|\le \sum^{k-1}_{j=0}\|p^{uk}_{j+1}-p^{uk}_j\|+\lm^k h_k\sum^{k-1}_{j=0}\|w^{uk}_j\|+2\sum^{k-1}_{j=0}|\xi^k_j|\to 0
$$
as $k\to\infty$ by the relationships in \eqref{tL} and \eqref{e:105a}.\vspace*{-0.05in}

To get a contradiction with our assumption on the violation of \eqref{e:83}, it remains by \eqref{e:97} to verify that $\xi^k_k\to 0$ as $k\to\infty$. To see this, observe that the convergence $p^{k}(T)\to 0,\;\lm^k\to 0$ implies by the first condition in \eqref{e:72} that $p^{xk}_k\to 0,\,p^{uk}_k\to 0$, and $\sum^m_{i=1}\eta^k_{ki}x^*_i\to 0$ as $k\to\infty$. Then it follows from the second condition in \eqref{e:72} that $\xi^k_k\ou^k_k\to 0$, which yields $\xi^k_k\to 0$ since $\ou^k_k\ne 0$ for large $k\in N$ due to $\ou^k_k\to\ou(T)$ and the assumptions on either $\ou(T)\ne 0$ or $\t<r$ that exclude vanishing $\ou^k_k$ by the constraints in \eqref{e:8}. Thus we arrive at a contradiction with \eqref{e:97} and so justify the nontriviality condition \eqref{e:83}.\\[1ex]
{\bf Step~8:} {\em Enhanced nontriviality conditions.} Our final step is to justify the stronger/enhanced nontriviality conditions in \eqref{enh1} and \eqref{enh2} under the interiority assumptions imposed therein provided that $0<\t<r$.\vspace*{-0.05in}

Consider first the left endpoint case \eqref{enh1} and suppose by contradiction that $(\lm,p(T))=0$ under the validity of the interiority condition in \eqref{enh1}. It follows from the latter that $1-\tau-\varepsilon_k<\|\ou^k(0)\|<1+\tau+\varepsilon_k$ for $i=1,\ldots,m$ and $k$ sufficiently large. Then we deduce from \eqref{gamma-k} and \eqref{e:70} that $\gamma^k_{0i}=0$ and $\xi^k_0=0$ for $i=1,\ldots,m$. Combining with \eqref{e:67} and the construction of $q^{uk}(\cdot)$ in Step~3 yields
$$
q^{uk}(0)=p^{uk}_0=p^{uk}_1-h_k\lm^kw^{uk}_0,\quad k\in\N.
$$
Since $p^{uk}_1\to 0$ as $k\to\infty$ by the above proof, we conclude that $q^u(0)=\lim_{k\to\infty}q^{uk}(0)=0$. This contradicts the nontriviality condition \eqref{e:83} and thus verifies \eqref{enh1}. The justification of \eqref{enh2} is similar, and therefore we complete the proof of the theorem. $\h$

Let us now specify the general necessary optimality conditions of Theorem~\ref{Th:12} to the important novel case of our consideration in this paper, where we have controls only in perturbations while $u$-controls in $(P^\t)$ are {\em fixed}. Such a setting is used in Section~5 for applications to the controlled crowd motion model. In this case each problem $(P)$ reduces to the following form $(\Tilde P)$:
\begin{equation*}
\mbox{minimize }\;\Tilde{J}[x,a]:=\varphi(x(T))+\int^T_0\ell\big(t,x(t),a(t),\dot{x}(t),\dot{a}(t)\big)dt
\end{equation*}
subject to the sweeping differential inclusion
\begin{equation}\label{e:5a}
 -\dot x(t)\in N\big(x(t)-\ou(t);C\big)+f\big(x(t),a(t)\big)\;\mbox{ a.e. }\;t\in[0,T],\quad x(0):=x_0\in C
\end{equation}
with the convex polyhedron $C$ in \eqref{e:6} and the implicit state constraints
\begin{equation*}
\big\la x^*_i,x(t)-\ou(t)\big\ra\le 0\;\mbox{ for all }\;t\in[0,T]\;\mbox{ and }\;i=1,\ldots,m,
\end{equation*}
which follow from \eqref{e:5a}. As above, we study problem $(\Tilde P)$ in the class of $(x(\cdot),a(\cdot))\in W^{1,2}([0,T];\R^{n+d})$. Observe that we do not need to consider in this case the $\t$-parametric version of $(\Tilde P)$.\vspace*{-0.05in}

The next result follows from the specification of Theorem~\ref{Th:12} and its proof in the case of $(\Tilde P)$ by taking into account the structures of the sweeping set $C(t)$ and the running cost $\ell$ therein.\vspace*{-0.1in}

\begin{corollary}{\bf(necessary conditions for sweeping optimal solutions with controlled perturbations).}\label{Th:15} Let $(\ox(\cdot),\oa(\cdot))$ be a given r.i.l.m.\ for $(\Tilde P)$ satisfying {\rm(H1)}, {\rm(H2)}, all the appropriate assumptions of Theorem~{\rm\ref{Th:11}}, and the assumptions on the running cost $\ell$ from Theorem~{\rm\ref{Th:12}} with $\ell_2=0$. Then there exist a number $\lm\ge 0$, subgradient functions $w(\cdot)=(w^x(\cdot),w^a(\cdot))\in L^2([0,T];\R^{n-d})$ and $v(\cdot)=(v^x(\cdot),v^a(\cdot))\in L^2(]0,T];\R^{n+d})$ well defined at $t=0$ and satisfying \eqref{subg}, an adjoint arc $p(\cdot)=(p^x(\cdot),p^a(\cdot))\in W^{1,2}([0,T];R^{n+d})$, and a measure $\gamma=(\gg_1,\ldots,\gg_n)\in C^*([0,T];\R^n)$ on $[0,T]$ such that we have conditions \eqref{e:76}, \eqref{e:74} with the functions $\eta_i\in L^2([0,T];\R_+)$ uniquely defined by representation \eqref{e:76} together with the following relationships for a.e.\ $t\in[0,T]$:
\begin{equation}\label{e:117}
\left\{\begin{array}{ll}
\dot{p}^x(t)=\lm w^x(t)+\nabla_x f\big(\ox(t),\oa(t)\big)^*\big(\lm v^x(t)-q^x(t)\big),\\
\dot{p}^a(t)=\lm w^a(t)+\nabla_a f\big(\ox(t),\oa(t)\big)^*\big(\lm v^x(t)-q^x(t)\big),
\end{array}\right.
\end{equation}
where the vector function $q=(q^x,q^a)\colon[0,T]\to\R^n\times\R^d$ is of bounded variation on $[0,T]$ satisfying
\begin{equation}\label{e:118}
q^a(t)\in\lm\partial\ell_3\big(t,\dot\oa(t)\big)\;\mbox{ for a.e.}\;t\in[0,T]\quad\mbox{and}
\end{equation}
\begin{equation}\label{e:119}
q(t):= p(t)-\int_{[t,T]}\big(d\gamma(s),0\big)\;\mbox{ on }\;[0,T]
\end{equation}
except at most a countable subset for its left-continuous representative. We also have the measure nonatomicity condition {\rm (a)} of Theorem~{\rm\ref{Th:12}} and the transversality relationships
\begin{equation}\label{e:121}
\left\{\begin{array}{ll}
-p^x(T)\in\lm\partial\varphi\big(\ox(T)\big)+\disp\sum_{i\in I(\ox(T)-\ou(T))}\eta_i(T)x^*_i\subset\partial\varphi\big(\ox(T)\big)+N\big((\ox(T)-\ou(T);C\big),\\
p^a(T)=0,\;\mbox{ and }\;q^a(0)=\lm v^a(0).
\end{array}\right.
\end{equation}
Finally, the enhanced nontriviality condition
\begin{equation}\label{e:122}
\lm+\|p(T)\|\ne 0
\end{equation}
holds provided that either $\la\ox(t),\ou(t)\ra\not=\|\ou(t)\|^2$ on $[0,T)$, or $\la x^*_i,x_0-\ou(0)\ra<0$ for all $i=1,\ldots,m$.
\end{corollary}\vspace*{-0.15in}
\begin{remark} {\bf (discussions on optimality conditions).}\label{disc} {\rm The results of Theorem~\ref{Th:12} and Corollary~\ref{Th:15} provide comprehensive necessary optimality conditions for a broad class of intermediate (between weak and strong with including the latter) local minimizers of state-constrained sweeping control problems concerning highly unbounded and non-Lipschitzian differential inclusions. Now we briefly discuss some remarkable features of the obtained results with their relationships to known results in this direction.\\[1ex]
{\bf (i)} As has been well recognized in standard optimal control theory for differential equations and Lipschitzian differential inclusions with state constraints, necessary optimality conditions for such problems may exhibit the {\em degeneration phenomenon} when they hold for every feasible solution with some nontrivial collection of dual variables. In particular, this could happen if the initial vector at $t=0$ belongs to the boundary of state constraints; see \cite{AS,v} for more discussions and references.
\vspace*{-0.05in}

It may also be the case for our problem $(P^\t)$ under the general nontriviality condition \eqref{e:83} when, e.g., $\t=0$ and the vector $x_0-\ou(0)$ lays on the boundary of the polyhedral set $C$. However, the degeneration phenomenon is surely excluded in $(P^\t)$ by the enhanced nontriviality in \eqref{enh1}, \eqref{enh2} and by the condition $\la\ox(t),\ou(t)\ra\not=\|\ou(t)\|^2$ on $[0,T)$ in $(\Tilde P)$ even in the case where either $x_0-\ou(0)$ or $\ox(T)-\ou(T)$ is a boundary point of the convex polyhedron $C$; see Examples~\ref{Ex:2}, \ref{Ex:3} and also Examples~\ref{cr2}, \ref{cr3} for the crowd motion model. As other examples demonstrate (see, in particular, Example~\ref{Ex:1}), even in the case of potential degeneracy as in \eqref{e:83} for $\t=0$ under the violation of the aforementioned conditions that rule out the degeneration phenomenon, the obtained results can be useful to determine optimal solutions.\\[1ex]
{\bf (ii)} Let us draw the reader's attention to some specific features of the new {\em transversality conditions} obtained in Theorem~\ref{Th:12} and Corollary~\ref{Th:15}. The transversality condition at the {\em left endpoint} in \eqref{e:81a} and \eqref{e:121} may look surprising at the first glance since the initial vector $x_0$ of the feasible sweeping trajectories $x(\cdot)$ is fixed. However, it is not the case for control functions $u(\cdot)$ and $a(\cdot)$, which are incorporated into the differential inclusion \eqref{e:20} and the cost functional \eqref{e:4} with their initial points being reflected in \eqref{e:81a}. The usage of the left transversality condition \eqref{e:81a} allows us to exclude in Example~\ref{Ex:1} the potential degeneration term $q^u(0)$ from the general nontriviality condition \eqref{e:83} and then to calculate an optimal solution to the sweeping control problem under consideration.\vspace*{-0.05in}

Observe that we get the two types of the transversality conditions for $p(T)$ at the {\em right endpoint} in \eqref{e:81} and \eqref{e:121}: one expressed directly via $\eta_i(T)$ and other given via the normal cone $N(\ox(T)-\ou(T);C)$ due to \eqref{e:81b}. While the second type of transversality is more expected, the first type is essentially more precise. Indeed, the normal cone transversality may potentially lead us to degeneration when $\ox(T)-\ou(T)$ lays at the boundary of $C$. On the other hand, degeneration is completely excluded in this case if we have $\eta_i(T)=0$ as $i\in I(\ox(T)-\ou(T))$ for the endpoint vectors $\eta_i(T)$, which may occur independently of their {\em a priori} location at $N(\ox(T)-\ou(T);C)$ due to the fact that the vectors $\eta_i(T)$ are {\em uniquely determined} by representation \eqref{e:76} of the term $-\dot{\ox}(T))-f(\ox(t),\oa(t))$ via the linearly independent generating vectors $x^*_i$. This is explicitly illustrated by Example~\ref{Ex:1}.\\[1ex]
{\bf (iii)} It has been largely understood in optimal control of differential equations and Lipschitzian differential inclusions that necessary optimality conditions for problems with inequality state constraints are described via {\em nonnegative Borel measures}. In the case of $(P^\t)$ we have both inequality and equality state constraints on $z(\cdot)$ given by \eqref{e:8} and \eqref{mixed} that are reflected in Theorem~\ref{Th:12} by the measure $\xi\in C^*([0,T];\R)$ and $\gg\in C^*([0,T];\R^m)$, respectively. In problem $(\Tilde P)$ we do not have state constraints for the $u$-components, and so only the measure $\gg$ appears in the optimality conditions of Corollary~\ref{Th:15}. But even in the latter case we {\em do not ensure the nonnegativity} of $\gg$ (see Examples~\ref{cr2} and \ref{cr3} for the controlled crowd motion model), which once more reveals a significant difference between the sweeping control problems governed in fact by {\em evolution/differential variational inequalities} from the conventional state-constrained control problem considered in the literature. On the other hand, all the examples presented in Sections~4 and 5 show that our results agree with those known for conventional models while indicating that the corresponding measures become nonzero at the points where optimal trajectories {\em hit the boundaries} of state constraints and {\em stay such} on these boundaries; see Examples~\ref{cr2} and \ref{cr3} to illustrate the latter phenomenon.\\[1ex]
{\bf (iv)} Finally, we compare the results derived in Theorem~\ref{Th:12} (and their consequences in Corollary~\ref{Th:15}) with the most recent necessary optimality conditions obtained in \cite[Theorem~6.1]{chhm2} for problem $({\bar P}^\t)$ as $0<\t<T$ of minimizing the cost functional
\begin{equation*}
\bar J[x,u,b]:=\varphi\big(x(T)\big)+\int_0^{T}\ell\big(t,x(t),u(t),b(t),\dot x(t),\dot u(t),\dot b(t)\big)dt
\end{equation*}
over absolutely continuous controls $u=(u_1,\ldots,u_m)\colon[0,T]\to\R^{mn}$,\;$b=(b_1,\ldots,b_m)\colon[0,T]\to\R^m$ and the corresponding absolutely continuous trajectories $x\colon[0,T]\to\R^n$ of the unperturbed sweeping inclusion
\begin{equation*}
-\dot x(t)\in N\big(x(t);C(t)\big)\;\mbox{ for a.e. }\;t\in[0,T],\quad x(0):=x_{0}\in C(0),
\end{equation*}
with the sweeping set $C(t)$ and the constraints on $u$-controls given by
\begin{equation*}
C(t):=\big\{x\in\R^n\big|\;\la u_i(t),x\ra\le b_i(t)\;\mbox{ for all }\;i=1,\ldots,m\big\},
\end{equation*}
\begin{equation*}
\|u_i(t)\|=1\;\mbox{ for all }\;t\in[0,T]\;\mbox{ and }\;i=1,\ldots,m.
\end{equation*}
We can see that problem $({\bar P}^\t)$ is different from $(P^\t)$ by the {\em absence of controlled perturbations} (which is of course the underlying feature of our problem $(P^\t)$ and its applications to the crowd motion model), the choice of $\t$, and a bit different class of feasible solutions. On the other hand, while ignoring these differences, problem $(P^\t)$ can be reduced to $\bP$ with no $u$-controls (they are replaced by the generating vectors $x^*_i$ of the polyhedron $C$) and with $b$-controls given in the form
$$
b_i(t):=\big\la x^*_i,u(t)\big\ra\;\mbox{ for all }\;t\in[0,T]\;\mbox{ and }\;i=1,\ldots,m
$$
via the $u$-controls in $(P^\t)$. However, problem $\bP$ obtained in this way from $(P^\t)$ in the absence of perturbations is not considered in \cite{chhm2}, since we do have the {\em pointwise constraints} on $u_i(t)$ in \eqref{e:8}, which are in fact a part of the state constraints on $z(t)$ in the setting of \eqref{e:20} under investigation, while there are {\em no any constraints} on $b_i(t)$ in \cite{chhm2}. Necessary optimality conditions for problems $\bP$ of this type (where $\t$ does not play any role since the $u$-controls are fixed) are specified in \cite[Theorem~6.3]{chhm2}. It is not hard to check that the results obtained therein are included in those established in Theorem~\ref{Th:12} for $(P^\t)$ in the case where both problems are the same. However, even in this (not so broad) case we obtain additional information in Theorem~\ref{Th:12} and Corollary~\ref{Th:15} in comparison with \cite{chhm2}. Let us list the main new ingredients of our results for $(P^\t)$ in the common setting with \cite[Theorem~6.3]{chhm2} and also in a similar (while different) setting of \cite[Theorem~6.1]{chhm2} for $\bP$ with $u$-control components, which can be incorporated therein by using the {\em more precise} discrete approximation technique developed in this paper:\vspace*{-0.05in}

$\bullet$ The new transversality conditions at the left endpoint; see remark (ii) above.\vspace*{-0.05in}

$\bullet$ Both types of transversality at the right endpoint discussed in remark (ii) are different and more convenient for applications in comparison with (6.10)--(6.12) in \cite{chhm2}. Observe that the latter ones are given implicitly as equations for $p^x,p^u,p^b$ at the local optimal solution $\oz(T)$.\vspace*{-0.05in}

$\bullet$ Our results are applied to the general case of the parameter $\t$ and its interrelation with another parameter $r$ in the $u$-control bounds in contrast to only the interior case of $\t\in(0,T)$  with $r=1$ in \cite{chhm2}.\vspace*{-0.05in}

$\bullet$ Our general nontriviality condition \eqref{e:83} contains only the $u$-component $q^u(0)$ in contrast to all the components of $q(0)$ in the corresponding condition $\lm+\|p(T)\|+\|q(0)\|\ne 0$ of \cite{chhm2}.\vspace*{-0.05in}

$\bullet$ Theorem~\ref{Th:12} and Corollary~\ref{Th:15} present more conditions that surely rule out the degeneracy phenomenon in comparison with the corresponding results of \cite[Theorems~6.1, 6.3]{chhm2}; see the discussion in remark (i). Note that the appearance of degeneracy is also excluded by the new transversality conditions as discussed in remark (ii) and illustrated by the examples below.\vspace*{-0.05in}

$\bullet$ The presence of controlled perturbations in $(P^\t)$ and $(\Tilde P)$ allows us to reveal new behavior phenomena for the measure $\gg$ responsible for the state constraints \eqref{mixed} in comparison with the settings of \cite{chhm2}, even in the absence of the measure $\xi$ responsible for the $u$-constraints in \eqref{e:8}; see remark (iii). In particular, Examples~\ref{cr2} and \ref{cr3} illustrate behavior of the measure $\gg$ in keeping the optimal trajectory on the boundary of state constraints in the crowd motion model.}
\end{remark}\vspace*{-0.3in}

\section{Numerical Examples}
\setcounter{equation}{0}\vspace*{-0.1in}

In this section we present three academic examples illustrating some characteristic features of the obtained necessary optimality conditions for problems $(P^\t)$ and $(\Tilde P)$ and their usefulness to determine optimal solutions and exclude nonoptimal ones in rather simple settings. More involved examples with our major applications to the crowd motion model in a corridor are given in Section~5.\vspace*{-0.1in}

\begin{example}\label{Ex:1} {\bf(optimal controls in both sweeping set and perturbations).} Consider problem $(P^\t)$ with any $0\le\t<1/2$ and the following data:
\begin{equation}\label{e:123}
n=m=d=T=1,\;x_0:=0,\;x^*_1:=1,\;f(x,a):=a,\;\varphi(x):=\dfrac{(x-1)^2}{2},\;\ell(t,x,u,a,\dot{x},\dot{u},\dot{a}):=\dfrac{1}{2}a^2.
\end{equation}
\begin{center}
\begin{tikzpicture}
\draw [thick] (1,-0.8) circle (0.1);
\draw node[below] at (1,-1) {0};
\draw [thick, red](0,-1)--(12,-1);
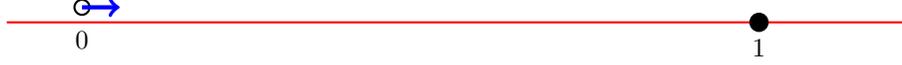
\draw node[below] at (6,-1.5) {{\bf Figure 1:} Direction of optimal control};
\draw [ultra thick, blue] [->] (1,-0.8)--(1.5,-0.8);
\draw [ultra thick, fill = black] (10,-1) circle (0.1);
\draw node[below] at (10,-1.1) {1};
\end{tikzpicture}
\end{center}
\vspace*{-0.1in}

In this case we have $C=\R_-$. The structure of the cost functional in \eqref{e:123} allows us to assume without loss of generality that $a$-controls are uniformly bounded, and thus $(P^\t)$ admits an optimal solution $(\ox(\cdot),\ou(\cdot),\oa(\cdot))\in W^{1,2}([0,1];\R^3)$ by \cite[Theorem~4.1]{cm1}. It is also easy to see that all the assumptions of Theorem~\ref{Th:12} are satisfied. Furthermore, it follows from the structure of $(P^\t)$ with $r=1/2$ in \eqref{e:8}  that $\ou(t)=1/2$ on $[\t,1-\t]$ and $\ou(t)\in[1/2-\t,1/2+\t]$ on $[0,\t)\cup(1-\t,1]$; see Figure~1. Supposing further that $\ox(t)\in{\rm int}(C+\ou(t))$ for any $t\in[0,1)$ and that $-\dot\ox(1)=f(\ox(1),\oa(1))$, we see that these assumptions are realized for the optimal solution found via the necessary optimality conditions of Theorem~\ref{Th:12}.\vspace*{-0.05in}

With taking into account that the second assumption in \eqref{assum-nontr} holds in our case, we get from Theorem~\ref{Th:12} the following relationships, where $\lm\ge 0$ and $\eta(\cdot)\in L^2([0,1];\R_+)$ being well defined at $t=1$:\\[1ex]
(1) $\quad$ $w(t)=\big(0,0,\oa(t)\big),\;v(t)=0$ for a.e.\ $t\in[0,1]$;\\[1ex]
(2) $\quad$ $\ox(t)<\ou(t)\Longrightarrow\eta(t)=0$ for a.e.\ $t\in[0,1)$;\\[1ex]
(3) $\quad$ $\eta(t)>0\Longrightarrow q^x(t)=0$ for a.e.\ $t\in[0,1]$ including $t=1$;\\[1ex]
(4) $\quad$ $-\dot{\ox}(t)=\eta(t)+f\big(\ox(t),\oa(t)\big)=\eta(t)+\oa(t)$ for a.e.\ $t\in[0,1]$ including $t=1$;\\[1ex]
(5) $\quad$ $\big(\dot{p}^x(t),\dot{p}^u(t),\dot{p}^a(t)\big)=\big(0,0,\lm\oa(t)-q^x(t)\big)$ for a.e.\ $t\in[0,1]$;\\[1ex]
(6) $\quad$ $q^u(t)=0,\;q^a(t)=0$ for a.e.\ $t\in[0,1]$;\\
(7) $\quad$ $\big(q^x(t),q^u(t),q^a(t)\big)=\big(p^x(t),p^u(t),p^a(t)\big)-\left(\disp{\int_{[t,1]}d\gamma,\int_{[t,1]}2d\xi-d\gamma},0\right)$ for a.e.\ $t\in[0,1]$;\\
(8) $\quad$ $-p^x(1)=\lm\big(\ox(1)-1\big)+\eta(1),\;p^a(1)=0$;\\[1ex]
(9) $\quad$ $p^u(1)\in\eta(1)+2\ou(1)N\big(\ou(1);[1/2-\t,1/2+\t]\big)$;\\[1ex]
(10) $\quad$ $\eta(1)\in N\big(\ox(1)-\ou(1);C\big)$;\\[1ex]
(11) $\quad$ $q^u(0)\in-2\ou(0)N\big(\ou(0);[1/2-\tau,1/2+\tau]\big)+D^*G\big(x_0-\ou(0),-\dot{\ox}(0)-\oa(0)\big)\big(-q^x(0)\big)$;\\[1ex]
(12) $\quad$ $\lm+|q^u(0)|+|p(1)|\ne 0$.

Since $\ox(t)-\ou(t)\in{\rm int}\,C$ for all $t\in[0,1)$, the coderivative of the mapping $G(\cdot)=N(\cdot;C)$ disappears in the left endpoint transversality condition (11). Furthermore, we have $\eta(1)=0$ by (3) due to the assumption $-\dot\ox(1)=f(\ox(1),\oa(1))$ imposed on the optimal solution. Hence condition (10) holds automatically, and we arrive at the updated transversality relationships:
\begin{equation}\label{e:123a}
\left\{\begin{array}{ll}
-p^x(1)=\lm\big(\ox(1)-1\big),\;p^u(1)\in 2\ou(1)N\big(\ou(1);[1/2-\t,1/2+\t]\big),\\\\
q^u(0)\in-2\ou(0)N\big(\ou(0);[1/2-\tau,1/2+\tau]\big).
\end{array}\right.
\end{equation}
It follows from (5)--(7) that $p^x(\cdot)$ is a constant function on $[0,1]$ and that
\begin{equation}\label{e:125}
\lm\oa(t)=q^x(t)=p^x(1)-\int_{[t,1]}d\gamma\;\mbox{ for a.e. }\;t\in[0,1].
\end{equation}\vspace*{-0.22in}

The next assertion that holds in any finite-dimensional space is a consequence of the measure nonatomicity condition (a) of Theorem~\ref{Th:12}, which is essential in this and other examples.\\[1ex]
{\bf Claim:} {\em Let $\la x^*,\ox(s)-\ou(s)\ra<0$ for all $s\in[t_1,t_2]$ with $t_1,t_2\in[0,T)$ and some vector $x^*\in\R^n$ under the validity of the measure nonatomicity condition $($a$)$ of Theorem~{\rm\ref{Th:12}} involving the vector $x^*$ and the measure $\gg$ therein. Then $\gamma([t_1,t_2])=0$ and $\gamma(\{s\})=0$ whenever $s\in[t_1,t_2]$. Thus we also have $\gamma((t_1,t_2))=\gamma([t_1,t_2))=\gamma((t_1,t_2])=0$.}

To verify this claim, pick any $s\in[t_1,t_2]$ with $\la x^*_1,\ox(t)-\ou(t)\ra<0$ and find by the imposed measure nonatomicity condition a neighborhood $V_s$ of $s$ in $[0,T]$ such that $\gamma(V)=0$ for all Borel subsets $V$ of $V_s$, and hence obviously $\gamma(\{s\})=0$. Since $[t_1,t_2]\subset\bigcup_{s\in[t_1,t_2]}V_s$ and $[t_1,t_2]$ is compact, there are finitely many points $s_1,\ldots,s_p\in[t_1,t_2]$ with $[t_1,t_2]\subset\bigcup_{i=1}^pV_{s_i}$. For each $i=1,\ldots,p-1$ take $\Tilde{s_i}\in V_{s_i}\cap V_{s_{i+1}}$ such that $[s_i,\Tilde{s_i}]\subset V_{s_i}$ and $[\Tilde{s_i},s_{i+1}]\subset V_{s_{i+1}}$, where $s_1=t_1$ and $s_p=t_2$ without loss of generality. Then the claim readily follows from the equalities
\begin{equation*}
\gamma([t_1,t_2])=\gamma\left(\bigcup_{i=1}^{p-1}[s_i,\Tilde{s_i})\cup[\Tilde{s_i},s_{i+1})\right)=\sum_{i=1}^{p-1}\Big(\gamma([s_i,\Tilde{s_i}))
+\gamma([\Tilde{s_i},s_{i+1}))\Big)=0.
\end{equation*}\vspace*{-0.15in}

Going back to our example, observe that the validity of $\ox(s)<\ou(s)$ for all $s\in[t,1)$ with $t\in[0,1)$ yields $\gamma([t,1])=\gamma(\{1\})$. Indeed, it follows from the above claim that for all large $k\in\N$ we get
\begin{equation*}
\begin{aligned}
\gamma([t,1])&=\gamma([t,1))+\gamma(\{1\})=\gamma\Big(\left[t,1-k^{-1}\right]\cup\bigcup_{n\ge k}\big(1-n^{-1},1-(n+1)^{-1}\big]\Big)+\gamma(\{1\})\\
&=\gamma\big(\big[t,1-k^{-1}\big]\big)+\sum_{n\ge k}\gamma\Big(\big(1-n^{-1},1-(n+1)^{-1}\big]\Big)+\gamma(\{1\})=\gamma(\{1\}).
\end{aligned}
\end{equation*}
This allows us to deduce from \eqref{e:125} that
\begin{equation}\label{at}
\lm\oa(t)=p^x(1)-\gamma(\{1\})\;\mbox{ for a.e. }\;t\in[0,1].
\end{equation}
\vspace*{-0.22in}

To proceed further, consider first the case where $1/2-\t<\ou(t)<1/2+\t$ for $t=0,1$. In this case we have $q^u(0)=p^u(1)=0$ by \eqref{e:123a} and so $\lm>0$, since the opposite would contradict the nontriviality condition (12) by taking (8) with $\eta(1)=0$ into account. It follows now from \eqref{at} that
$\oa(\cdot)$ must be a constant function, $\oa(\cdot)\equiv\vartheta$ on $[0,1]$, due to its continuity. Then (2) and (4) ensure that
$$
\ox(t)=x_0+\int_0^t\dot{\ox}(s)ds=-\int_0^t\vartheta ds=-t\vartheta\;\mbox{ for all }\;t\in[0,1].
$$
Consequently, the cost functional in our problem is calculated by
$$
J[\ox,\ou,\oa]=\dfrac{(-\vartheta-1)^2}{2}+\dfrac{{\vartheta}^2}{2}=\vartheta^2+\vartheta+\dfrac{1}{2}
$$
and clearly achieves its absolute minimum at $\vartheta=-1/2$. Thus in this case we arrive by the necessary optimality conditions of Theorem~\ref{Th:12} at the (local) optimal solution
$$
\ox(t)=t/2,\;\;\oa(t)=-1/2\;\mbox{ on }\;[0,1],\;\;\ou(t)=1/2\;\mbox{ on }\;[\tau,1-\tau],\;\mbox{ and }\;\ou(t)\in(1/2-\t,1/2+\t)\;\mbox{ on }\;[0,\t)\cup(1-\t,1],
$$
which satisfies all the preliminary assumptions imposed above.\vspace*{-0.05in}

In the case where $\ou(t)\in\{1/2-\t,1/2+\t\}$ as $t=0,1$ we get from \eqref{e:123a} that
\begin{equation*}
\left\{\begin{array}{ll}
p^u(1)\le0,\;q^u(0)\ge 0&\mbox{if}\quad\ou(0)=1/2-\t,\;\ou(1)=1/2-\t,\\
p^u(1)\ge0,\;q^u(0)\ge 0&\mbox{if}\quad\ou(0)=1/2-\t,\;\ou(1)=1/2+\t,\\
p^u(1)\le0,\;q^u(0)\le 0&\mbox{if}\quad\ou(0)=1+\t,\;\ou(1)=1-\t,\\
p^u(1)\ge0,\;q^u(0)\le 0&\mbox{if}\quad\ou(0)=1/2+\t,\;\ou(1)=1/2+\t,
\end{array}\right.
\end{equation*}
which does not provide sufficient information to conclude that $p^u(1)=q^u(0)=0$ and thus $\lm>0$. If the latter holds, we can proceed similarly to the interior case and find local minimizers as above. However, the case of $\lm=0$ remains open from the viewpoint of Theorem~\ref{Th:12}.\vspace*{-0.05in}

Observe finally that by fixing the $u$-control component as $\ou(\cdot)=1/2$ on $[0,1]$ we reduce problem $(P^\t)$ of this example to the type of $(\Tilde P)$. Then the necessary conditions of Corollary~\ref{Th:15} and the arguments above allow us to calculate, by taking into account the existence result of \cite[Theorem~4.1]{cm1}, the unique global optimal solution $(\ox(t),\oa(t))=(t/2,-1/2)$ for all $t\in[0,1]$.
\end{example}\vspace*{-0.12in}

The next example concerns problem $(\Tilde P)$ with controlled perturbations and nonsmooth data while illustrating the usage of necessary optimality conditions of Corollary~\ref{Th:15} to determine optimal solutions and designate nonoptimal ones under different parameter values.\vspace*{-0.1in}

\begin{example}{\bf (nonsmooth problems with controlled perturbations).}\label{Ex:2}
Consider problem $(\Tilde P)$ with $n,m,d,T,x_0,x^*_1,f(x,a)$ as in \eqref{e:123}, fixed $\ou(t)=r$ on $[0,1]$, and the cost functions $\varphi(x):=(x-1)^2$,
\begin{equation}\label{ex4.2}
\ell(t,x,u,a,\dot{x},\dot{u},\dot{a}):=(a+2t)^2+\alpha|\dot{a}+4t-1|\;\mbox{ for }\;\alpha\ge 0.
\end{equation}\vspace*{-0.25in}

Let us first examine the case of the parameters $r=1$ and $\al=0$. The structure of the cost functional in this case suggests a natural candidate for the optimal solution $(\ox(t),\oa(t))=(t^2,-2t)$ on $[0,1]$. Observe that $\ox(t)<\ou(t)$ for all $t\in[0,1)$ and that $\ox(t)\cdot\ou(t)\not=|\ou(t)|^2 $ for all $t\in[0,1)$. Applying now the necessary optimality conditions of Corollary~\ref{Th:15} gives us the following relationships with a number $\lm\ge 0$ and a function $\eta(\cdot)\in L^2([0,1];\R_+)$ well defined at $t=1$:\\[1ex]
(1) $\quad$ $w(t)=\big(0,2(\oa(t)+2t)\big),\;v(t)=(0,0)$ for a.e.\ $t\in[0,1]$;\\[1ex]
(2) $\quad$ $\ox(t)<\ou(t)\Longrightarrow\eta(t)=0$ for a.e.\ $t\in[0,1]$ including $t=1$;\\[1ex]
(3) $\quad$ $\eta(t)>0\Longrightarrow q^x(t)=0$ for a.e.\ $t\in[0,1]$;\\[1ex]
(4) $\quad$ $-\dot{\ox}(t)=\eta(t)+f\big(\ox(t),\oa(t)\big)=\eta(t)+\oa(t)$ for a.e.\ $t\in[0,1]$;\\[1ex]
(5) $\quad$ $\big(\dot{p}^x(t),\dot{p}^a(t)\big)=\big(0,2\lm(\oa(t)+2t)-q^x(t)\big)$ for a.e.\ $t\in[0,1]$;\\[1ex]
(6) $\quad$ $q^a(t)=0$ for a.e.\ $t\in[0,1]$;\\
(7) $\quad$ $\big(q^x(t),q^a(t)\big)=\big(p^x(t),p^a(t)\big)-\left(\disp{\int_{[t,1]}d\gamma},0\right)$ for a.e.\ $t\in[0,1]$;\\
(8) $\quad$ $-p^x(1)\in 2\lm\big(\ox(1)-1\big)+\eta(1),\;p^a(1)=0;$\\[1ex]
(9) $\quad$ $\eta(1)\in N\big(\ox(1)-\ou(1);C\big)$;\\[1ex]
(10)$\quad$ $\lm+|p^x(1)|\not=0$.

Combining the relationships in (5)--(7) gives us the equation
\begin{equation}\label{e:126}
2\lm\big(\oa(t)+2t\big)=q^x(t)=p^x(1)-\gamma([t,1])\;\mbox{ for a.e. }\;t\in[0,1].
\end{equation}
Letting $\lm>0$ and taking into account that $x(t)-u(t)<0$ for all $t\in[0,1)$, we get by the arguments in Example~\ref{Ex:1} that
$\gamma([t,1])=\gg(\{1\})$. It implies that $\oa(t)+2t$ reduces to a constant $\vt$ a.e.\ on $[0,1]$, which ensures that $\oa(t)=-2t+\vartheta$ for all $t\in [0,1]$ due to the continuity of $\oa(\cdot)$. It follows from (4) that
$$
\ox(t)=\int_0^t\dot{\ox}(s)ds =-\int_0^t\eta(s)ds-\int_0^t(-2s+\vartheta)ds=t^2-t\vartheta,\quad t\in[0,1].
$$
Thus $\ox(1)=1-\vartheta$, which gives the value of 2$\vartheta^2$ to the cost functional with the minimal value achieved at $\vt=0$. This confirms via Corollary~\ref{Th:15} the optimality of the solution $(\ox(t),\oa(t))=(t^2,-2t)$ for the above choice of the parameters $(r,\al)=(1,0)$ in the problem under consideration. Note that the other conditions in (1)--(10) besides those used above hold automatically for $(\ox(\cdot),\oa(\cdot))$ with $p^x(1)=\eta(1)=0$.\vspace*{-0.05in}

Consider now this problem $(\Tilde P)$ with the parameter values $r=2$ and $\al>0$; the latter generates nonsmoothness in \eqref{ex4.2}. Let us check by the feasible solution $(\ox(t),\oa(t))=(t^2,-2t)$ on $[0,1]$ is not locally optimal anymore for $(\Tilde P)$ by using the necessary optimality conditions of Corollary~\ref{Th:15} listed above with the replacement of (6) by the subdifferential inclusion \eqref{e:118} in the nonsmooth case of \eqref{ex4.2}. It follows from (7) that $q^a(t)=p^a(t)$ for a.e.\ $t\in[0,1]$. Furthermore, \eqref{e:118} tells us in this case that $q^a(t)=-\alpha\lm$ for a.e.\ $0\le t<3/4$ and $q^a(t)=\alpha\lm$ for a.e.\ $3/4<t\le 1$. Thus we have
\begin{equation*}
p^a(t)=\left\{\begin{array}{ll}
-\alpha\lm&\mbox{for all }\;0\le t<3/4,\\
\alpha\lm&\mbox{for all }\;3/4<t\le 1,
\end{array}\right.
\end{equation*}
which yields by the continuity of $p^a(\cdot)$ on $[0,1]$ that $\al\lm=0$ and hence $\lm=0$. Then it follows from (8) and (2) with $\ox(1)=1<\ou(1)=2$ in that $p^x(1)=0$. This contradicts the nontriviality condition (10) and hence justifies that the given pair $(\ox(\cdot),\oa(\cdot))$ fails to be an optimal solution to $(\Tilde P)$ with $r=2$ and $\al>0$.
\end{example}\vspace*{-0.12in}

The the next example addresses a two-dimensional perturbed sweeping process and demonstrates the possibility to determine optimal solutions by using the necessary optimality conditions of Corollary~\ref{Th:15}.\vspace*{-0.1in}

\begin{example}{\bf (two-dimensional sweeping process with controlled perturbations.)}\label{Ex:3} Consider problem $(\Tilde{P})$ with the following initial data:
$$
n=m=d=2,\;T=1,\;x_0:=(0,-1),\;x^*_1:=(1,0),\;x^*_2:=(0,1),\;f(x,a):=a,\;\varphi(x):=0,
$$
and $\ell(t,x,u,a,\dot{x},\dot{u},\dot{a}):=(\|\dot{x}\|^2+\|a\|^2)/2$. Given $\ou(\cdot)=(1,0)$ on $[0,1]$, apply the necessary optimality conditions of Corollary~\ref{Th:15} to determine (local) optimal solutions $\oa(\cdot)=(\oa_1(\cdot),\oa_2(\cdot))$ and $\ox(\cdot)=(\ox_1(\cdot),\ox_2(\cdot))$ to this problem. We seek for solutions to $(\Tilde P)$ such that
\begin{equation}\label{ex4.3}
\la x^*_i,\ox(t)-\ou(t)\ra<0\;\mbox{ for all }\;t\in[0,1),\;i=1,2,\;\mbox{ and }\;\ox(1)-\ou(1)\in{\rm bd}\,C.
\end{equation}
and show that \eqref{ex4.3} holds for $(\ox(\cdot),\oa(\cdot))$ found below by using the necessary optimality conditions of Corollary~\ref{Th:15}. In the case of $(\Tilde P)$ under consideration these conditions look as follows, where $\lm\ge 0$ and $\eta(\cdot)=\big(\eta_1(\cdot,\eta_2(\cdot)\big)\in L^2([0,1];\R_+)$ well defined at $t=1$:\\[1ex]
(1) $\quad$ $w(t)=\big(0,0,\oa(t)\big),\;v(t)=\big(\dot{\ox}(t),0,0\big)$ for a.e.\ $t\in[0,1]$;\\[1ex]
(2) $\quad$ $\la x^*_i,\ox(t)-\ou(t)\ra<0\Longrightarrow\eta_i(t)=0$ for $i=1,2$ and a.e.\ $t\in[0,1)$;\\[1ex]
(3) $\quad$ $\eta_i(t)\Longrightarrow\la x^*_i,\lm\dot{\ox}(t)-q^x(t)\ra=0$ for $i=1,2$ and  a.e.\ $t\in[0,1]$;\\[1ex]
(4) $\quad$ $-\dot{\ox}(t)=\big(-\dot{\ox}_1(t),-\dot{\ox}_2(t)\big)=(\eta_1(t),\eta_2(t))+\big(\oa_1(t),\oa_2(t)\big)$ for a.e.\ $t\in[0, 1]$;\\[1ex]
(5) $\quad$ $\big(\dot{p}^x(t),\dot{p}^a(t)\big)=\lm\big(0,\oa(t)\big)+\big(0,(\lm\dot{\ox}_1(t)-q^x_2(t),\lm\dot{\ox}_2(t)-q^x_2(t))\big)$ for a.e.\ $t\in[0, 1]$;\\[1ex]
(6) $\quad$ $q^a(t)=0$ for a.e.\ $t\in[0,1]$;\\[1ex]
(7) $\quad$ $q^x(t)=p^x(t)-\gamma([t,1]),\quad q^a(t)=p^a(t)$ for a.e.\ $t\in[0,1]$;\\[1ex]
(8) $\quad$ $-p^x(1)=\big(\eta_1(1),\eta_2(1)\big)\in N\big(\ox(1)-\ou(1);C\big)$;\\[1ex]
(9) $\quad$ $\lm+\|p^x(1)\|\ne 0$.

Employing the first condition in \eqref{ex4.3} together with (2) and (4), we obtain that $\dot{\ox}(t)=-\oa(t)$ for a.e.\ $t\in[0,1]$. It also follows from (5)--(7) that
\begin{equation}\label{e:128}
\lm\oa(t)=\lm\dot{\ox}(t)-q^x(t),\;\mbox{ i.e., }\;2\lm\oa(t)=-q^x(t)\;\mbox{ for a.e. }\;t\in[0,1]
\end{equation}
Using (5) again tells us that $p^x(\cdot)$ is constant on $[0,1]$, i.e., $p^x(t)\equiv p^x(1)$. This allows us to deduce that
\begin{equation*}
q^x(t)=p^x(1)-\gamma([t,1])=p^x(1)-\gamma(\{1\})\;\mbox{ for a.e. }\;t\in[0,1]
\end{equation*}
by using the measure nonatomicity condition (a) from Theorem~\ref{Th:12} for the measure $\gg$ and repeating the arguments of Example~\ref{Ex:1}. This shows by \eqref{e:128} and the control continuity that $\oa(\cdot)$ is a constant on $[0,1]$ provided that $\lm\ne 0$; otherwise, we do not have enough information to proceed. Putting $\oa(t)\equiv(\vartheta_1,\vartheta_2)$ for all $t\in[0,1]$ gives us $\ox(t)=(-\vartheta_1t,-1-\vartheta_2t)$ for all $t\in[0,1]$. Thus $\ox(1)=(-\vartheta_1,-1-\vartheta_2)$, and by the second condition in \eqref{ex4.3} we have the following two possibilities:\vspace*{-0.05in}

{\bf (a):} $\ox_1(1)=1$. Then $\vartheta_1=-1$ and the cost functional reduces is $J[\ox,\oa]=1+\vartheta_2^2$. It obviously achieves its absolute minimum value $\bar J=1$ at the point $\vartheta_2=0$.

{\bf (b):} $\ox_2(1)=0$. Then $\vartheta_2=-1$, and the minimum cost is $\bar J=1$ achieved at $\vartheta_1=0$.\vspace*{-0.05in}

As a result, we arrive at are two feasible solutions giving the same optimal value to the cost functionals:
$$
\ox(t)=(t,-1),\;\oa(t)=(-1,0)\quad\mbox{and}\quad\ox(t)=(0,t-1),\;\oa(t)=(0,-1).
$$
Figure~2 provides some illustration of the sweeping motion in this case, where the red lines indicate the boundary points at which the corresponding sweeping trajectories hit the state constraints.
\end{example}
\begin{center}
\begin{tikzpicture}
\draw [thick] (7,-2) circle (0.1);
\draw node[below] at (7,-2.2) {$x_0=(0,-1)$};
\draw node[above] at (7,-1) {(0,0)};
\draw node[right] at (8,-2) {(1,-1)};
\draw [thick, red](3,-1)--(8,-1);
\draw [thick, red](8,-1)--(8,-4);
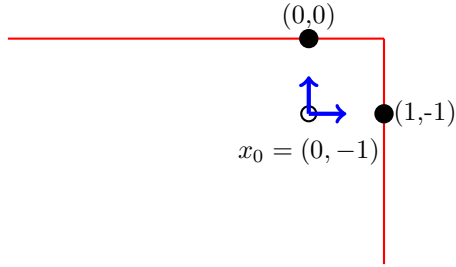
\draw node[below] at (6,-4) {{\bf Figure~2:} Two-dimensional motion.};
\draw [ultra thick, blue] [->] (7,-2)--(7.5,-2);
\draw [ultra thick, blue] [->] (7,-2)--(7,-1.5);
\draw [ultra thick, fill = black] (7,-1) circle (0.1);
\draw [ultra thick, fill = black] (8,-2) circle (0.1);
\end{tikzpicture}
\end{center}
\vspace*{-0.15in}

\section{Controlled Crowd Motion Model in a Corridor}
\setcounter{equation}{0}\vspace*{-0.1in}

This section is devoted to the formulation and solution of an {\em optimal control problem} concerning the so-called {\em crowd motion model in a corridor}. We refer the reader to \cite{mv0,mv,ve} for describing of the dynamics in such and related crowd motion models as a {\em sweeping process} with the corresponding mathematical theory, numerical simulations, and various applications. However, neither these papers nor other publications contain, to the best of our knowledge, control and/or optimization versions of crowd motion modeling, which is of our main interest here. We follow the terminology and notation of \cite{mv0,mv,ve}.\vspace*{-0.03in}

Our main goal is to demonstrate that the necessary optimality conditions obtained in Corollary~\ref{Th:15} allow us to develop an effective procedure to determine optimal solutions in the general setting under consideration with finitely many participants and then explicitly solve the problem in some particular situations involving two and three participants. Furthermore, in this way we reveal certain specific features of the obtained necessary optimality conditions for problems with state constraints.

The crowd motion model of \cite{mv0,mv,ve} is designed to deal with local interactions between participants in order to describe the dynamics of pedestrian traffic. This model rests on the following postulates:\vspace*{-0.05in}

$\bullet$ A spontaneous velocity corresponding to the velocity that each participant would like to have in the absence of others is defined first.\vspace*{-0.05in}

$\bullet$ The actual velocity is then calculated as the projection of the spontaneous velocity onto the set of admissible velocities, i.e., such velocities that do not violate certain nonoverlapping constraints.\vspace*{-0.03in}

In what follows we consider $n$ participants $(n\ge 2$) identified with rigid disks of the same radius $R$ in a corridor as depicted in Figure~3.
\begin{center}
\begin{tikzpicture}
\draw [thick, red](0,0)--(12,0);
\draw [thick] (1,-0.5) circle (0.5);
\draw node[below] at (1.1,-0.5) {$x_1$};
\draw [thick] (2.5,-0.5) circle (0.5);
\draw node[below] at (2.6,-0.5) {$x_2$};
\draw [thick] (4,-0.5) circle (0.5);
\draw node[below] at (4.1,-0.5) {...};
\draw [thick] (5.2,-0.5) circle (0.5);
\draw node[below] at (5.3,-0.5) {$x_i$};
\draw [thick] (6.2,-0.5) circle (0.5);
\draw node[below] at (6.3,-0.5) {$x_{i+1}$};
\draw [thick] (7.2,-0.5) circle (0.5);
\draw node[below] at (7.3,-0.5) {...};
\draw [thick] (8.2,-0.5) circle (0.5);
\draw node[below] at (8.3,-0.5) {$x_{N}$};
\draw [thick, red](0,-1) --(12,-1);
\draw node[below] at (6,-1.5) {{\bf Figure~3:} Crowd motion model in a corridor} ;
\draw [ultra thick, blue] [->] (1,-0.5) -- (1.5,-0.5);
\draw [ultra thick, blue] [->] (2.5,-0.5) -- (3.2,-0.5);
\draw [ultra thick, blue] [->] (4,-0.5) -- (4.6,-0.5);
\draw [ultra thick, blue] [->] (5.2,-0.5) -- (5.7,-0.5);
\draw [ultra thick, blue] [->] (6.2,-0.5) -- (6.5,-0.5);
\draw [ultra thick, blue] [->] (7.2,-0.5) -- (7.5,-0.5);
\draw [ultra thick, blue] [->] (8.2,-0.5) -- (8.5,-0.5);
\draw [ultra thick, fill = orange] (10,-1) rectangle (11,0);
\draw node[below] at (10.5,-1) {Exit};
\end{tikzpicture}
\end{center}\vspace*{-0.1in}

In that case, since the participants are not likely to leap across each other, it is natural to restrict the set of feasible configurations to one of its connected components ({\em nonoverlapping condition}):
\begin{equation}\label{e:129}
Q_0=\big\{x=(x_1,\ldots,x_n)\in\R^n,\;x_{i+1}-x_i\ge 2R\big\}.
\end{equation}
Assuming that the participants exhibit the same behavior, their {\em spontaneous velocity} can be written as
$$
U(x)=\big(U_0(x_1),\ldots,U_0(x_n)\big)\;\mbox{ for }\;x\in Q_0,
$$
where $Q_0$ is taken from \eqref{e:129}. Observe that the nonoverlapping constraint in \eqref{e:129} does not allow the participants to move with their spontaneous velocity, and the distance between two participants in contact can only increase. To reflect this situation, the set of {\em feasible velocities}
$$
C_x:=\big\{v=(v_1,\ldots,v_n)\in\R^n\left|\right.x_{i+1}-x_i=2R\Longrightarrow v_{i+1}\ge v_i\;\mbox{ for all }\;i=1,\ldots,n-1\big\},
$$
and then describe the {\em actual velocity field} is the feasible field via the Euclidean projection of $U(x)$ to $C_x$:
$$
\dot{x}(t)=\Pi\big(U(x);C_x\big)\;\mbox{ for a.e. }\;t\in[0,T],\quad x(0)=x_0\in Q_0,
$$
where $T>0$ is a fixed duration of the process and $x_0$ indicates the starting position of the participants. Using the orthogonal decomposition via the sum of mutually polar cone as in \cite{mv,ve}, we get
\begin{equation*}
U(x)\in N_x+\dot{x}(t)\;\mbox{ for a.e. }\;t\in[0,T],\quad x(0)=x_0,
\end{equation*}
where $N_x$ stands for the normal cone to $Q_0$ at $x$ and can be described in this case as the polar
$$
N_x=C^*_x:=\big\{w\in\R^n\left|\right.\la w,v\ra\le 0\;\mbox{ for all }\;v\in C_x\big\},\quad x\in Q_0.
$$\vspace*{-0.22in}

Let us now rewrite this model in the form used in our problem $(\Tilde P)$ without control parameters so far. Given the orths $(e_1,\ldots,e_n)\in\R^n$, specify the polyhedral set $C$ by
\begin{equation}\label{e:131}
C:=\big\{x\in\R^n\left|\right.\la x^*_i,x\ra\le 0,\;i=1,\ldots,n-1\big\}\;\mbox{ with }\;x^*_i:=e_i-e_{i+1},\quad i=1,\ldots,n-1.
\end{equation}
Since all the participants exhibit the same behavior and want to reach the exit by the shortest path, their spontaneous velocities can be represented as
\begin{equation*}
U(x)=\big(U_0(x_1),\ldots,U_0(x_n)\big)\;\mbox{ with }\;U_0(x)=-s\nabla D(x),
\end{equation*}
where $D(x)$ stands for the distance between the position $x=(x_1,\ldots,x_n)\in Q_0$ and the exit, and where the scalar $s\ge 0$ denotes the speed. Since $x\ne 0$ and hence $\|\nabla D(x)\|=1$, we have $s=\|U_0(x)\|$. By taking this into account as well as the aforementioned postulate that, in the absence of other participants, each participant tends to remain his/her spontaneous velocity until reaching the exit, the (uncontrolled) perturbations in this model are described by
\begin{equation*}
f(x)=-(s_1,\ldots,s_n)\in\R^n\;\mbox{ for all }\;x=(x_1,\ldots,x_n)\in Q_0,
\end{equation*}
where $s_i$ denotes the speed of the participant $i=1,\ldots,n$. However, if participant $i$ is in contact with participant $i+1$ in the sense that $x_{i+1}(t)-x_i(t)=2R$, then both of them tend to adjust their velocities in order to keep the distance at least $2R$ with the participant in contact. To control the actual speed of all the participants in the presence of the nonoverlapping condition \eqref{e:129}, we suggest to involve
control functions $a(\cdot)=(a_1(\cdot),\ldots,a_n(\cdot))$ into perturbations as follows:
\begin{equation}\label{cr-pert}
f\big(x(t),a(t)\big)=\big(s_1a_1(t),\ldots,s_n a_n(t)\big),\quad t\in[0,T].
\end{equation}
In order to represent this controlled crowd motion model in the form of $(\Tilde P)$, define recurrently the vector function $\ou=(\ou_1,\ldots,\ou_n)\colon[0,T]\to\R^n$, which is constant in our case, by
\begin{equation}\label{cr-u}
\ou_{i+1}(t)-\ou_i(t)=2R\;\mbox{ with }\;\ou_1(t)=\alpha\;\mbox{ and }\;\|\ou(t)\|=r,\quad i=1,\ldots,n-1,
\end{equation}
where $r=r(\al)$ is an increasing function of $\al$ with the value of $\al$ specified later. Note that the nonoverlapping condition \eqref{e:129} can be written now, due to \eqref{e:131} and \eqref{cr-u}, via the state constraints
\begin{equation}\label{cr-state}
x(t)-\ou(t)\in C\;\mbox{ for all }\;t\in[0,T],
\end{equation}
where the points $t\in[0,T]$ with $x_{i+1}(t)-x_i(t)=2R$ are exactly those at which the motion $x(t)-\ou(t)$ hits the polyhedral constraint set $C$.\vspace*{-0.05in}

The constructions above allow us to present the {\em controlled crowd motion dynamics} as
\begin{equation}\label{e:135}
\left\{\begin{array}{lcl}
-\dot{x}(t)\in N\big(x(t);C(t)\big)+f\big(x(t),a(t)\big)\;\mbox{ for a.e. }\;t\in[0,T],\\
C(t):=C+\ou(t),\;\|\ou(t)\|=r\;\mbox{ on }\;[0,T],\quad x(0)=x_0\in C(0),
\end{array}\right.
\end{equation}
with $C$, $f$, and $\ou$ taken from \eqref{e:131}, \eqref{cr-pert}, and \eqref{cr-u}, respectively. Recall that the state constraints \eqref{cr-state} are implicitly present in \eqref{e:135} due to definition \eqref{nor} of the normal cone to convex sets.

To optimize dynamics \eqref{e:135} by using controls $a(\cdot)$, we introduce the {\em cost functional}
\begin{eqnarray}\label{cr-cost}
\mbox{minimize }\;J[x,a]:=\dfrac{1}{2}\Big(\|x(T)\|^2+\int^T_0\|a(t)\|^2dt\Big)
\end{eqnarray}
the meaning of which is to {\em minimize the distance} of all the participants to the exit at the origin together with the {\em energy} of feasible controls $a(\cdot)$. Having now the formulated optimal control problem for the crowd motion model in the form of $(\Tilde P)$, we can apply to solving this problem the necessary optimality conditions for the sweeping process with controlled perturbations derived in Corollary~\ref{Th:15}.\vspace*{-0.03in}

It is easy to see that all the assumptions of Corollary~\ref{Th:15} are satisfied for problem \eqref{e:135}, \eqref{cr-cost}. To make sure that the nontriviality condition holds in the enhanced/nondegenerate form \eqref{e:122}, we select the parameter $\al$ in \eqref{cr-u} so large that
$$
r=r(\al)>l=\|x_0\|+e^{2MT}2MT\big(1+\|x_0\|\big)
$$
where the number $l>0$ is calculated in \eqref{et-es} for the constant control $u(\cdot)$. As mentioned in \eqref{assum-nontr} of Theorem~\ref{Th:12}, this condition with $\t=0$ yields the validity of the second condition therein, which ensures in turn the fulfillment of the enhanced nontriviality \eqref{e:122} in Corollary~\ref{Th:15}.\vspace*{-0.05in}

Applying now the necessary optimality conditions of Corollary~\ref{Th:15} gives us the following, where $\lm\ge 0$ and $\eta_i(\cdot)\in L^2([0,T];\R_+)$ well defined at $t=T$:\\[1ex]
(1) $\quad$ $w(t)=\big(0,\oa(t)\big),\;v(t)=(0,0)$ for a.e.\ $t\in[0,T]$;\\
(2) $\quad$ $-\dot{\ox}(t)=\disp{\sum^{n-1}_{i=1}}\eta_i(t)x^*_i+(s_1\oa_1(t),\ldots,s_n\oa_n(t))$ for a.e.\ $t\in[0,T];$\\
(3) $\quad$ $\ox_{i+1}(t)-\ox_i(t)>2R\Longrightarrow\eta_i(t)=0$ for all $i=1,\ldots,n-1$ and a.e.\ $t\in[0,T]$;\\[1ex]
(4) $\quad$ $\eta_i(t)>0\Longrightarrow q^x_i(t)=q^x_{i+1}(t)$ for all $i=1,\ldots,n-1$ and a.e.\ $t\in[0,T]$;\\[1ex]
(5) $\quad$ $\dot{p}(t)=\big(0,\lm\oa(t)-(s_1q^x_1(t),\ldots,s_nq^x_n(t))\big)$ for a.e.\ $t\in[0,T]$;\\[1ex]
(6) $\quad$ $q^x(t)=p^x(t)-\gamma([t,T])$ for a.e.\ $t\in[0,T]$;\\[1ex]
(7) $\quad$ $q^a(t)=p^a(t)=0$ for a.e.\ $t\in[0,T]$;\\[1ex]
(8) $\quad$ $-p^x(T)=\lm\ox(T)+\sum_{i\in I(\ox(T)-\ou(T))}\eta_i(T)$;\\[1ex]
(9) $\quad$ $\sum_{i\in I(\ox(T)-\ou(T))}\eta_i(T)\in N\big(\ox(T)-\ou(T);C)$;\\[1ex]
(10) $\quad$ $p^a(T)=0$;\\[1ex]
(11) $\quad$ $\lm+\|p^x(T)\|\not=0$.

As discussed above, the situation where $\ox_{i+1}(t_1)-\ox_i(t_1)=2R$ for some $t_1\in[0,T]$ pushes participants $i$ and $i+1$ to adjust their speeds in order to keep the distance at least $2R$ with the one in contact. It is natural to suppose that both participants $i$ and $i+1$ maintain their new constant velocities after the time $t=t_1$ until either reaching someone ahead or the end of the process at time $t=T$. Hence the velocities of all the participants are piecewise constant on $[0,T]$ in this setting.\vspace*{-0.05in}

Observe that the differential relation in (2) can be read as
\begin{equation}\label{e:135a}
\left\{\begin{array}{ll}
-\dot{\ox}_1(t)=\eta_1(t)+s_1\oa_1(t),\\
-\dot{\ox}_i(t)=\eta_i(t)-\eta_{i-1}(t)+s_i\oa_i(t),\quad i=2,\ldots,n-1,\\
-\dot{\ox}_n(t)=-\eta_{n-1}(t)+s_n\oa_n(t)
\end{array}\right.
\end{equation}
for a.e.\ $t\in[0,T]$. Next we clarify the sense of the implications in (3). If participant~1 is far away from participant~2 in the sense that $\ox_2(t)-\ox_1(t)>2R$ for some time $t\in[0,T]$, then his/her actual velocity and the spontaneous velocity are the same meaning that $-\dot{\ox}_1(t)=s_1\oa_1(t)$. Likewise we have the same situation for the last participant $n$. However, it is not the case for two adjacent participants between the first and last ones because they must rely on the participants ahead and behind them.\vspace*{-0.05in}

Further, it follows from condition (5) that we have
\begin{equation}\label{cr-lm}
\lm\oa_i(t)=s_iq^x_i(t)\mbox{ for a.e. }\;t\in[0,T]\;\mbox{ and all }\;i=1,\ldots,n.
\end{equation}
If $\eta_i(t)>0$ for some $i\in\{1,\ldots,n-1\}$ and $t\in[0,T]$, we deduce from (4) and \eqref{cr-lm} by taking into account the continuity of $\oa_i(\cdot)$ on $[0,T]$ that
\begin{equation}\label{e:136}
s_{i+1}\oa_i(t)=s_i\oa_{i+1}(t)\;\mbox{ for all }\;t\in[0,T]\;\mbox{ and }\;i=1,\ldots,n-1
\end{equation}
provided that $\lm>0$ (say $\lm=1$); otherwise, we do not have enough information to proceed. Since the velocities $s_i$ are constant in \eqref{e:137}, it is to suppose by \eqref{e:136} that the functions $\oa_i(\cdot)$ are constant $\oa_i$ on $[0,T]$ for all $i=1,\ldots,n$ and thus optimal controls among such functions. Using this and the Newton-Leibniz formula in \eqref{e:135a} gives us the trajectory representations for all $t\in[0,T]$:
\begin{equation}\label{e:137}
\left\{\begin{array}{ll}
\ox_1(t)=x_{01}-\disp\int^t_0\eta_1(s)ds-ts_1\oa_1,\\
\ox_i(t)=x_{0i}+\disp\int^t_0\big[\eta_{i-1}(s)-\eta_i(s)\big]ds-ts_i\oa_i\;\mbox{ for }\;i=2,\ldots,n-1,\\
x_n(t)=x_{0n}+\disp\int^t_0\eta_{n-1}(s)ds-ts_n\oa_n,
\end{array}\right.
\end{equation}
where $(x_{01},\ldots,x_{0n})$ are the components of the starting point $x_0\in\R^n$ in \eqref{e:135}.

Prior to developing an effective procedure to find optimal solutions to the controlled crowd motion model by using the obtained optimality conditions in the general case above, we consider the following example for two participants that shows how to explicitly solve the problem in such settings.\vspace*{-0.1in}

\begin{example} {\bf(solving the crowd motion control problem with two participants).}\label{cr2} Specify the data of \eqref{e:135}, \eqref{cr-cost} as follows:
$n=2,\;T=6,\;s_1=6,\;s_2=3,\;x_{01}=-60,\;x_{02}=-48,\;R=3$. Then the equations in \eqref{e:137} reduce for all $t\in[0,6]$ to
\begin{equation}\label{e:139}
\ox_1(t)=-60-\disp{\int^t_0}\eta(s)ds-6t\oa_1,\quad\ox_2(t)=-48+\disp{\int^t_0}\eta(s)ds-3t\oa_2.
\end{equation}
Let $t_1\in[0,6]$ be the first time when $\ox_2(t_1)-\ox_1(t_1)=2R=6$; see Figure~4.
\begin{center}
\begin{tikzpicture}
\draw [thick, red](0,0)--(12,0);
\draw [thick] (1,-0.5) circle (0.5);
\draw node[below] at (1.1,-0.5) {$x_1$};
\draw [thick] (3.5,-0.5) circle (0.5);
\draw node[below] at (3.6,-0.5) {$x_2$};
\draw [thick, red](0,-1) --(12,-1);
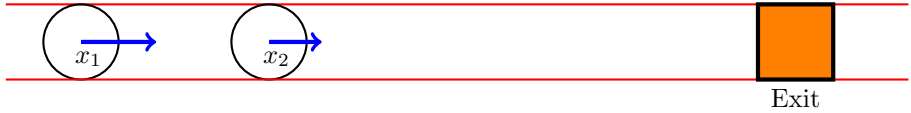
\draw node[below] at (6,-1.5) {{\bf Figure~4:} Two participants out of contact for $t<t_1$.};
\draw [ultra thick, blue] [->] (1,-0.5) -- (2,-0.5);
\draw [ultra thick, blue] [->] (3.5,-0.5) -- (4.2,-0.5);
\draw [ultra thick, fill = orange] (10,-1) rectangle (11,0);
\draw node[below] at (10.5,-1) {Exit};
\end{tikzpicture}
\end{center}\vspace*{-0.1in}

Hence for $t<t_1$ we have $x_2(t)-x_1(t)>2R=6$, and so $\eta(t)=0$ by (3). Note that at the point $t=t_1$ the motion $\ox(t)-\ou(t)$ hits the state constraint set $C$ in \eqref{cr-state} and thus is reflected by a nonzero measure $\gg$ in (6). However, we can proceed by an easier way in our particular setting. Indeed, subtracting the first equation in \eqref{e:139} from the second one gives us the relationship
\begin{equation}\label{oa}
12-3t_1(\oa_2-2\oa_1)=6\;\mbox{ and thus }\;6\oa_1-3\oa_2+1\le 0.
\end{equation}
Suppose without loss of generality that both functions $\eta(t)$  and $\dot\ox(t_1)$ are well defined at $t=t_1$. Then we get from \eqref{e:135a} and \eqref{e:139} in this case the expressions \begin{equation*}
\dot{\ox}_1(t_1)=-\eta(t_1)-6\oa_1\;\mbox{ and }\;\dot{\ox}_2(t_1)=\eta(t_1)-3\oa_2
\end{equation*}
with $\dot{\ox}_1(t_1)\le\dot{\ox}_2(t_1)$, which imply in turn that
\begin{equation}\label{eta}
-2\eta(t_1)-6\oa_1+3\oa_2\le 0.
\end{equation}
It follows from \eqref{oa} and \eqref{eta} that $\eta(t_1)\ge 1/2$. Furthermore, we deduce from \eqref{e:136} with the chosen speed values $s_1,s_2$ that the constant controls
$\oa_1,\oa_2$ are related by
$$
\oa_1=\dfrac{s_1}{s_2}\oa_2=2\oa_2.
$$
Having in hand the relationships above, let us now calculate an optimal solution to the problem under consideration by imposing the requirement that both participants maintain their new constant velocities until the end of the process at $t=T$, i.e., satisfying the condition $\dot{\ox}(t)=\dot{\ox}(t_1)$ for all $t\in[t_1,6]$. Since $\oa(\cdot)$ is constant on $[0,6]$ and $\dot{\ox}(\cdot)$ is constant on the intervals $[0,t_1)$ and $[t_1,6]$, the vector function $\eta(\cdot)$ is constant on $[0,t_1)$ and $[t_1,6]$ while admitting by \eqref{e:137} the representation
\begin{equation*}
\eta(t)=\left\{\begin{array}{ll}
\eta(0)&\mbox{a.e. }\;t\in[0,t_1)\;\mbox{ including }\;t=0,\\
\eta(t_1)&\mbox{a.e. }\;t\in[t_1,6]\;\mbox{ including }\;t=t_1.
\end{array}\right.
\end{equation*}
In particular, $\eta(t)=\eta(t_1)>0$ a.e.\ on $[t_1,6]$, and thus we get from (3) that $\ox_2(t)-\ox_1(t)=2R=6$ for all $t\in[t_1,6]$, i.e., the optimal motion stays on the boundary of state constraints \eqref{cr-state} on the whole interval $[t_1,6]$ meaning that the two participants of the model are in contact on this interval; see Figure~5.
\begin{center}
\begin{tikzpicture}
\draw [thick, red](0,0)--(12,0);
\draw [thick] (3.5,-0.5) circle (0.5);
\draw node[below] at (3.6,-0.5) {$x_1$};
\draw [thick] (4.5,-0.5) circle (0.5);
\draw node[below] at (4.6,-0.5) {$x_2$};
\draw [thick, red](0,-1) --(12,-1);
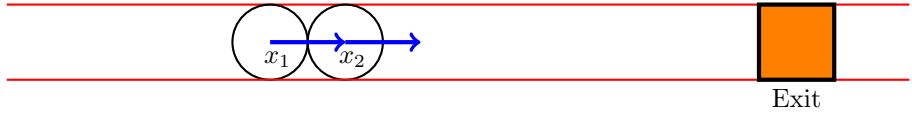
\draw node[below] at (6,-1.5) {{\bf Figure 5:} Two participants in contact for $t\ge t_1$.};
\draw [ultra thick, blue] [->] (3.5,-0.5) -- (4.5,-0.5);
\draw [ultra thick, blue] [->] (4.5,-0.5) -- (5.5,-0.5);
\draw [ultra thick, fill = orange] (10,-1) rectangle (11,0);
\draw node[below] at (10.5,-1) {Exit};
\end{tikzpicture}
\end{center}\vspace*{-0.1in}

Combining this with the the subtraction of the first equation from the second one in \eqref{e:139} gives us
$$
(t-t_1)\big[2\eta(t_1)+6\oa_1-3\oa_2\big]=0\;\mbox{ for all}\;t\in[t_1,6],
$$
which in turn implies that $2\eta(t_1)+6\oa_1-3\oa_2=0$. Remembering that $\oa_1=2\oa_2$, we calculate the value of $\eta(\cdot)$ at the hitting point $t=t_1$ by $\eta(t_1)=-\frac{9}{2}\oa_2=-\frac{9}{4}\oa_1$. Note also that $\dot{\ox}_2(t_1)=\dot{\ox}_1(t_1)$ in our case. Based on these calculations, we can express the value of cost functional \eqref{cr-cost} for this example at $(\ox,\oa)$ as
$$
J[\ox,\oa]=\frac{1}{2}\Big[\big(45\oa_2+57\big)^2+\big(45\oa_2+51\big)^2\Big]+15\oa^2_2.
$$
Minimizing this function of $\oa_2$ subject to the constraint $\oa_2\le-\frac{1}{9}$ that comes from the second expression in \eqref{oa} gives us the optimal control value $\oa_2=-\frac{4860}{4080}\approx-1.1911$. Accordingly the formulas obtained above allows us to calculate all the other ingredients of the optimal solution with the corresponding values of dual variables in the necessary optimality condition. It gives us, in particular, that
$$
\gg\big([t,6]\big)\approx(-1.56,3.76)\;\mbox{ for }\;0.56=t_1\le t\le 6,
$$
which reflects the fact that the optimal sweeping motion hits the boundary of the state constraints at $t_1=0.56$ and stays there till the end of the process at $T=6$. It is worth mentioning that the obtained nonzero measure $\gg$ has the opposite signs of its components on $[t_1,6]$, which is different from the standard optimal control problems with inequality state constraints.
\end{example}\vspace*{-0.1in}

Now we come back to the {\em general case} of the controlled crowd model in a corridor with $n\ge 3$ participants. Following the approach employed in Example~\ref{cr2}, we develop an {\em effective procedure} to determine an optimal control from the obtained necessary optimality conditions and then fully implement by a numerical example for the case where $n=3$. \vspace*{-0.05in}

Recall our postulate that any two adjacent participants $i$ and $i+1$ that come to be in contact at some point $t\in[0,T]$ (i.e., $x_{i+1}(t)-x_i(t)=2R$) have the same velocity therein, change their velocities at the contact point, and maintain their new constant velocities until reaching the participant ahead or until the end of the process at $t=T$.  This yields that the function $\eta(\cdot)$ in the conditions above is piecewise constant on $[0,T]$. Suppose for simplicity that $\eta_0(t)=\eta_n(t)=0$ on $[0,T]$ and then rewrite \eqref{e:137} as
$$
\ox_i(t)=x_{0i}+\disp{\int^t_0}\big[\eta_{i-1}(s)-\eta_i(s)\big]ds-ts_i\oa_i\;\mbox{ for }\;i=1,\ldots,n.
$$
Fix $i\in\{1,\ldots,n-1\}$, and let $t_i$ be the first time when $\ox_{i+1}(t_i)-\ox_i(t_i)=2R$; see Figure~6.
\begin{center}
\begin{tikzpicture}
\draw [thick, red](0,0)--(12,0);
\draw [thick] (1,-0.5) circle (0.5);
\draw node[below] at (1.1,-0.5) {$x_1$};
\draw [thick] (3.5,-0.5) circle (0.5);
\draw node[below] at (3.6,-0.5) {$x_2$};
\draw [thick] (4.5,-0.5) circle (0.5);
\draw node[below] at (4.6,-0.5) {$x_3$};
\draw [thick, red](0,-1) --(12,-1);
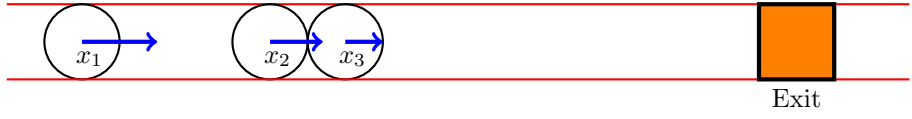
\draw node[below] at (6,-1.5) {{\bf Figure~6:} Out of contact situation for two adjacent participants when $t<t_1$};
\draw [ultra thick, blue] [->] (1,-0.5) -- (2,-0.5);
\draw [ultra thick, blue] [->] (3.5,-0.5) -- (4.2,-0.5);
\draw [ultra thick, blue] [->] (4.5,-0.5) -- (5,-0.5);
\draw [ultra thick, fill = orange] (10,-1) rectangle (11,0);
\draw node[below] at (10.5,-1) {Exit};
\end{tikzpicture}
\end{center}\vspace*{-0.1in}

For each such index $i$ consider the numbers
\begin{equation}\label{e:143a}
\vTh^i:=\min\big\{t_j\big|\;t_j>t_i,\;j=1,\ldots,n-1\big\},\quad\vTh_i:=\max\big\{t_j\big|\;t_j<t_i,\;j=1,\ldots,n-1\big\}
\end{equation}
and observe the following relationships for the optimal crowd motion on the intervals $[0,t_i)$ and $\in[t_i,\vTh^i)$:\vspace*{-0.05in}

$\bullet$ If $t\in[0,t_i)$, we have $\eta_i(\cdot)=0$ on this interval by (3). This gives us
\begin{equation*}
\ox_i(t)=x_{0i}+\disp{\int^t_0}\eta_{i-1}(s)ds-ts_i\oa_i,\quad
\ox_{i+1}(t)=x_{0(i+1)}-\disp{\int^t_0}\eta_{i+1}(s)ds-ts_{i+1}\oa_{i+1}\;\mbox{ for }\;t\in[0,t_i).
\end{equation*}\vspace*{-0.1in}

$\bullet$ If $t\in[t_i,\vTh^i)$ with $\vTh^i$ from \eqref{e:143a}, we have on this interval that
\begin{equation*}
\left\{\begin{array}{ll}
\ox_i(t)=x_{0i}+\disp\int^{t_i}_0\eta_{i-1}(s)ds+(t-t_i)\big[\eta_{i-1}(t_i)-\eta_i(t_i)\big]-ts_i\oa_i,\\
x_{i+1}(t)=x_{0(i+1)}-\disp\int^{t_i}_0\eta_{i+1}(s)ds +(t-t_i)\big[\eta_i(t_i)-\eta_{i+1}(t_i)\big]-ts_{i+1}\oa_{i+1}.
\end{array}\right.
\end{equation*}
In what follows we suppose without loss of generality that the functions $\dot\ox(\cdot)$ are well defined at $t_i$ while the functions $\eta(\cdot)$ are well defined at $t_i$ and $\vt_i$. Since at the contact time $t=t_i$ the distance between the two participants $i$ and $i+1$ is exactly $2R$ (see Figure~7), we have the following relationships:
\begin{center}
\begin{tikzpicture}
\draw [thick, red](0,0)--(12,0);
\draw [thick] (3.5,-0.5) circle (0.5);
\draw node[below] at (3.6,-0.5) {$x_1$};
\draw [thick] (4.5,-0.5) circle (0.5);
\draw node[below] at (4.6,-0.5) {$x_2$};
\draw [thick] (5.5,-0.5) circle (0.5);
\draw node[below] at (5.6,-0.5) {$x_3$};
\draw [thick, red](0,-1) --(12,-1);
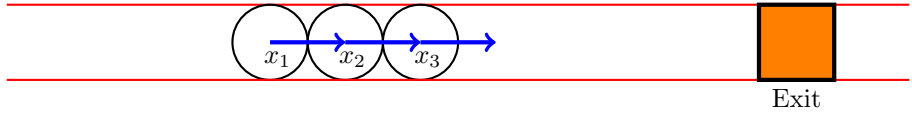
\draw node[below] at (6,-1.5) {{\bf Fig 7} All the participants in contact for $t\ge t_1$};
\draw [ultra thick, blue] [->] (3.5,-0.5) -- (4.5,-0.5);
\draw [ultra thick, blue] [->] (4.5,-0.5) -- (5.5,-0.5);
\draw [ultra thick, blue] [->] (5.5,-0.5) -- (6.5,-0.5);
\draw [ultra thick, fill = orange] (10,-1) rectangle (11,0);
\draw node[below] at (10.5,-1) {Exit};
\end{tikzpicture}
\end{center}\vspace*{-0.1in}
\begin{equation*}
\begin{aligned}
&2R=\ox_{i+1}(t_i)-\ox_i(t_i)=x_{0,(i+1)}-x_{0i}-\int^{t_i}_0\big[\eta_{i+1}(s)+\eta_{i-1}(s)\big]ds-t_i\big(s_{i+1}\oa_{i+1}-s_i\oa_i\big)\\
&=x_{0(i+1)}-x_{i0}-\int^{\vTh_i}_0\big[\eta_{i+1}(s)+\eta_{i-1}(s)\big]ds-(t_i-\vTh_i)\big[\eta_{i+1}(\vTh_i)+\eta_{i-1}(\vTh_i)\big],
-t_i(s_{i+1}\oa_{i+1}-s_ia_i),
\end{aligned}
\end{equation*}
where $\vTh_i$ is defined in \eqref{e:143a} being dependent of $t_i$. Then we can find $t_i\le T$ from the equation
\begin{equation}\label{e:145}
t_i=\dfrac{x_{0(i+1)}-x_{0i}-2R+\vTh_i\big[\eta_{i+1}(\vTh_i)+\eta_{i-1}(\vTh_i)\big]-\disp{\int^{\vTh_i}_0}\big[\eta_{i+1}(s)+\eta_{i-1}(s)\big]ds}{\eta_{i+1}(\vTh_i)
+\eta_{i-1}(\vTh_i)+s_{i+1}\oa_{i+1}-s_i\oa_i}
\end{equation}
provided that $x_{0(i+1)}-x_{0i}>2R$. In the case where $x_{0(i+1)}-x_{0i}=2R$ we put $t_i=0$. Our postulate tells us that $\dot{\ox}_{i+1}(t_i)=\dot{\ox}_{i}(t_i)$, which implies therefore that
\begin{equation}\label{e:147}
2\eta_i(t_i)=\eta_{i+1}(t_i)+\eta_{i-1}(t_i)+s_{i+1}\oa_{i+1}-s_i\oa_i.
\end{equation}
If $\eta_i(t_i)>0$, we get from the above that \eqref{e:136} holds, while the remaining case where $\eta_i(t_i)=0$ can be treated via \eqref{e:147}. The cost functional \eqref{cr-cost} can be expressed in this way as a function of $(\oa_1,\ldots,\oa_n)$ and $\eta_i(t_j)$ for $i=0,\ldots,n$ and $j=1,\ldots,n-1$. Consequently the optimal control problem under consideration reduces to the finite-dimensional optimization of this cost subject to inequality \eqref{e:145} and equality \eqref{e:147} constraints. To furnish these operations step-by-step, we proceed as follows:\\[1ex]
{\bf Step~1:} Determine which participants are in contact at the initial time, i.e., for which $i\in\{1,\ldots,n\}$ we have $x_{0(i+1)}-x_{0i}=2R$. If this occurs only for $i=n$, there is nothing to do. If it is the case of some $i\in\{0,\ldots,n-1\}$, we put $t_i:=0$ and observe that participants $i$ and $i+1$ have the same velocities while being away by $2R$ from each other.\\[1ex]
{\bf Step~2:} If $x_{0(i+1)}-x_{0i}=2R$ for all $i=0,\ldots,n-1$, express $t_i$ as a function of $\oa_i$ and $\eta_i$ by solving equation \eqref{e:145} for $t_i$ with $\vTh_i$ taken from in \eqref{e:143a}.\\[1ex]
{\bf Step~3:} Find relations between $\eta_i$ and $\oa_i$ from \eqref{e:136} and \eqref{e:147}, respectively, and substitute them into the cost function \eqref{cr-cost} for the subsequent optimization with respect to $\oa_i$.

We now demonstrate how this procedure work in the case where $n=3$ in the crowd motion model.\vspace*{-0.1in}

\begin{example} {\bf (solving the crowd motion control problem with three participants).}\label{cr3} Consider the optimal control problem in \eqref{e:135}, \eqref{cr-cost} with the following initial data:
$$
n=3,\;s_1=6,\;s_2=3,\;s_3=2,\;x_{01}=-60,\;x_{02}=-48,\;x_{03}=-42,\;T=6,\;R=3.
$$
By using the procedure outlined above, we first get $x_{02}-x_{01}=12>6=2R$ and $x_{03}-x_{02}=6=2R$. Then it is obvious that $t_2=0$, $t_1$ is determined by \eqref{e:145} as
$$
t_1=\dfrac{6}{\eta_2(0)+3\oa_2-6\oa_1}\le 6,
$$
and thus $\vartheta_1=t_2=0$. It is easy to see that in this example we have
\begin{equation*}
\dot{\ox}_1(t)=-6\oa_1,\;\dot{\ox}_2(t)=-\eta_2(0)-3\oa_2,\;\dot{\ox}_3(t)=\eta_2(0)-2\oa_3,
\end{equation*}
\begin{equation*}
\ox_1(t)=-60-6\oa_1,\;\ox_2(t)=-48-t\eta_2(0)-3t\oa_2,\;\ox_3(t)=-42+t\eta_2(0)-2t\oa_3
\end{equation*}
for $0\le t<t_1$, while for $t\in[t_1,6]$ the corresponding formulas are:
\begin{equation*}
\dot{\ox}_1(t)=-\eta_1(t_1)-6\oa_1,\;\dot{\ox}_2(t)=-\eta_2(t_1)+\eta_1(t_1)-3\oa_2,\;
\dot{\ox}_3(t)=\eta_2(t_1)-2\oa_3,
\end{equation*}
\begin{equation*}
\left\{\begin{array}{ll}
\ox_1(t)=-60-(t-t_1)\eta_1(t_1)-6t\oa_1,\;\ox_2(t)=-48-t\eta_2(0)+(t-t_1)\big(\eta_1(t_1)-\eta_2(t_1)\big)-3t\oa_2,\\
\ox_3(t)=-42+t\eta_2(0)+(t-t_1)\eta_2(t_1)-2t\oa_3.
\end{array}\right.
\end{equation*}
It follows directly from \eqref{e:147} the following relationships for $\eta(\cdot)$:
\begin{equation*}
2\eta_1(t_1)=\eta_2(t_1)+3\oa_2-6\oa_1,\;2\eta_2(0)=2\oa_3-3\oa_2,\;
2\eta_2(t_1)=\eta_1(t_1)+2\oa_3-3\oa_2.
\end{equation*}
Denoting for convenience $x:=a_2,\;y:=a_3,\;z:=a_1$ and taking into account that $\oa_1=2\oa_2$ by \eqref{e:136} due to $\eta_1(t_1)>0$, we rewrite these expressions and the above formula for $t_1$ as
\begin{equation}\label{e:151}
t_1=\dfrac{6}{-(21/2)x+y},\;\eta_1(t_1)=-8x+(13/6)y,\;
\eta_2(0)=-(3/2)x+y,\;\eta_2(t_1)=-3x+(4/3)y.
\end{equation}\vspace*{-0.1in}

Let us split the situation into the following two cases:\\[1ex]
{\bf Case~1: $\eta_2(0)>0$}. In this case we have $x=\frac{3}{2}y$, and thus \eqref{e:151} gives us the calculations:
\begin{equation*}
t_1=-(24/59y)\le 6,\;\eta_1(t_1)=-(59/6)y,\;\eta_2(0)=-(5/4)y,\;\eta_2(t_1)=-(37/6)y.
\end{equation*}
As a result, we have the expressions for the terminal points of the optimal trajectories
\begin{equation*}
\ox_1(6)=-49y-56,\;\ox_2(6)=-49y-50,\;\ox_3(6)=-49y-44
\end{equation*}
and the corresponding representation of the cost function
$$
J=\frac{1}{2}\Big[(49y+56)^2+(49y+50)^2+(49y+44)^2\Big]+36.75y^2.
$$
Minimizing this quadratic function over the constraint $y\le-\frac{4}{59}$ lead us to the optimal point $y=-\frac{7350}{7276.5}\approx-1.01$ and the corresponding values of the optimal control $\oa(t)=(-3.03,-1.52,-1.01)$, which creates the optimal contact time $t_1=0.40$ and the optimal crowd motion dynamics
\begin{equation*}
\big(\ox_1(t),\ox_2(t),\ox_3(t)\big)=\left\{
\begin{array}{ll}
(18.18t-60,\;3.28t-48,\;3.28t-42)&\mbox{for}\;t\in[0,t_1),\\
(8.25t-56,\;8.25t-50,\;8.25t-44)&\mbox{for}\;t\in[t_1,6].
\end{array}\right.
\end{equation*}
Note also that $\gamma([t,6])=(-2.92,\,4.71,\,1.24)$ when $t\in[t_1,6]$ with $\lm=1$ as considered above.\\[1ex]
{\bf Case~2:} $\eta_2(0)=0$. Then we can deduce from \eqref{e:151} that
\begin{equation*}
t_1=-y^{-1}\le 6,\;\eta_1(t_1)=-(19/6)y,\;\eta_2(t_1)=-(2/3)y,\;\eta_2(0)=-(3/2)x+y=0.
\end{equation*}
Hence $\eta_2(t_1)>0$, which implies by \eqref{e:136} that $2\oa_2=3\oa_3$ and so $x=\frac{3}{2}y$. Combining the latter with the above relation $x=\frac{2}{3}y$ tells us that $x=y=0$ This contradicts the constraint $y<0$ and thus rules out the situation in case. Overall, the calculations in Case~1 completely solve the crowd motion optimal control problem in this example by using the optimality conditions  established in Corollary~\ref{Th:15}.
\end{example}
\vspace*{-0.1in}

{\bf Acknowledgements.} The authors are very grateful to Juliette Venel for helpful discussions on the crowd motion model in a corridor and its applications.

\vspace*{-0.2in}

\begin{thebibliography}{99}

\bibitem{ao} L. Adam and J. V. Outrata. On optimal control of a sweeping process coupled with an ordinary differential equation. {\em Discrete
Contin. Dyn. Syst.--Ser. B}, 19:2709--2738, 2014.

\bibitem{AS} A. V. Arutyunov and S. M. Aseev. Investigation of the degeneracy phenomenon of the maximum principle for optimal control problems with
state constraints, {\em SIAM J. Control Optim.}, 35: 930--952, 1997.

\bibitem{ak} A. V. Arutyunov and D. Yu. Karamzin. Nondegenerate necessary optimality conditions for the optimal control problems with
equality type state constraints. To appear in {\em J. Global Optim.}, DOI:10.1007/s10898-015-0272-9, 2015.

\bibitem{bk} M. Brokate and P. Krej\v{c}\'\i. Optimal control of ODE systems involving a rate independent variational inequality.
{\em Discrete Contin. Dyn. Syst.--Ser. B}, 18:331--348, 2013.

\bibitem{cm1} T. H. Cao and B. S. Mordukhovich. Optimal control of a perturbed sweeping process via discrete approximations. Submitted for publication,
available in http://arxiv.org/abs/..., 2015.

\bibitem{clsw} F. H. Clarke, Yu. S Ledyaev, R. J. Stern and P. R. Wolenski. {\em Nonsmooth Analysis and Control Theory}. Springer, 1998.

\bibitem{chhm1} G. Colombo, R. Henrion, N. D. Hoang and B. S. Mordukhovich. Optimal control of the sweeping process. {\em Dyn. Contin. Discrete
Impuls. Syst.--Ser. B}, 19:117--159, 2012.

\bibitem{chhm2} G. Colombo, R. Henrion, N. D. Hoang and B. S. Mordukhovich. Optimal control of the sweeping process over polyhedral controlled
sets. To appear in {\em J. Diff. Eqs.}, http://dx.doi.org/10.1016/j.jde.2015.10.039, 2015.

\bibitem{CT} G. Colombo and L. Thibault. Prox-regular sets and applications. In Y. Gao and D. Motreanu, editors, {\em Handbook of Nonconvex Analysis},
pages 99--182. International Press, 2010.

\bibitem{dfm} T. Donchev, E. Farkhi and B. S. Mordukhovich. Discrete approximations, relaxation, and optimization of one-sided
Lipschitzian differential inclusions in Hilbert spaces. {\em J. Diff. Eqs.}, 243:301-328, 2007.

\bibitem{et} J. F. Edmond and L. Thibault. Relaxation of an optimal control problem involving a perturbed sweeping process. {\em Math. Program.},
104:347--373, 2005.

\bibitem{hmn} R. Henrion, B. S. Mordukhovich and N. M. Nam. Second-order analysis of polyhedral systems in finite and infinite dimensions
with applications to robust stability of variational inequalities. {\em SIAM J. Optim.}, 20:2199--2227, 2010.

\bibitem{Kr} P. Kre\v{c}\'\i. Evolution variational inequalities and multidimensional hysteresis operators. In P. Drabek, P. Kre\v{c}\'\i $\;$ and
P. Takac, editors, {\em Nonlinear Differential Equations}, Res. Notes Math. 404, pages 47--110. Chapman \& Hall, CRC, 1999.

\bibitem{KMM} M. Kunze and M. D. P. Monteiro Marques. An introduction to Moreau's sweeping process. In B. Brogliato, editor, {\em Impacts in
Mechanical Systems}, Lecture Notes in Phys. 551, pages 1--60. Springer, 2000.

\bibitem{mv0} B. Maury and J. Venel. A mathematical framework for a crowd motion model. {\em C. R. Acad. Sci. Paris Ser. I}, 346:1245--1250, 2008.

\bibitem{mv} B. Maury and J. Venel. Handling of contacts in crowd motion simulations. In C. Appert-Rolland et al., editors, {\em Traffic and Granular Flow '07},
pages 171--180. Springer, 2009.

\bibitem{mi} A. Mielke. Evolution of rate-independent systems. Evolutionary equations. Vol. II. In C. M. Dafermos and E. Feireis, editors,
{\em Handbook of Differential Equations}. Elsevier, 2005.

\bibitem{m95} B. S. Mordukhovich. Discrete approximations and refined Euler-Lagrange conditions for differential inclusions.
{\em SIAM J. Control Optim.}, 33:882--915, 1995.

\bibitem{m-book1} B. S. Mordukhovich. {\em Variational Analysis and Generalized Differentiation, I: Basic Theory}. Springer, 2006.

\bibitem{m-book2} B. S. Mordukhovich. {\em Variational Analysis and Generalized Differentiation, II: Applications}. Springer, 2006.

\bibitem{mor_frict} J. J. Moreau. On unilateral constraints, friction and plasticity. In G. Capriz and G. Stampacchia, editors, {\em New
Variational Techniques in Mathematical Physics}, Proceedings of C.I.M.E.\ Summer Schools, pages 173--322. Cremonese, 1974.

\bibitem{rw} R. T. Rockafellar and R. J-B. Wets. {\em Variational Analysis}. Springer, 2004.

\bibitem{smb} A. H. Siddiqi, P. Manchanda and M. Brokate. On some recent developments concerning Moreau's sweeping process. In A. H. Siddiqi and
 M. Ko\v{c}vara, editors, {\em Trends in Industrial and Applied Mathematics}, pages 339--354. Kluwer, 2002.

\bibitem{ve} J. Venel. A numerical scheme for a class of sweeping process. {\em Numerische Mathematik}, 118:451-484, 2011.

\bibitem{v} R. B. Vinter. {\em Optimal Control}, Birkha\"user, 2000.
\end{thebibliography}
\end{document}